\documentclass[aps,pre,twocolumn,groupedaddress,showpacs]{revtex4}

\usepackage{amsmath}
\usepackage{amsxtra}
\usepackage{amstext}
\usepackage{amssymb}

\usepackage{siunitx}

\usepackage{physics}

\usepackage{hyperref}
\hypersetup{
    colorlinks,%
    citecolor=black,%
    filecolor=black,%
    linkcolor=blue, %midnight-blue,%
    urlcolor=black,
    pagebackref=true,
}

\usepackage{graphicx}
\graphicspath{{figures/}}	
 
\newcommand{\eq}[1]{Eq.~(\ref{#1})}
\newcommand{\dde}[1]{DDE~(\ref{#1})}

\begin{document}

%\preprint{}

\title{Dynamics of a time-delayed relay system}

% Place the author information here.  Please hand-code the contact
% information and notecalls; do *not* use \footnote commands.  Let the
% author contact information appear immediately below the author names
% as shown.  We would also prefer that you don't change the type-size
% settings shown here.


\author{Lucas Illing}
\email[]{illing@reed.edu}
\affiliation{Department of Physics, Reed College, Portland, Oregon, 97202, USA}

\author{Pierce Ryan}
\affiliation{School of Mathematical Sciences, University College Cork, Cork T12 XF62, Ireland}

\author{Andreas Amann}
\affiliation{School of Mathematical Sciences, University College Cork, Cork T12 XF62, Ireland}

%\thanks{}
%\date{December 2017}

%\date{\today}
\begin{abstract}

We study the dynamics of a piecewise-linear second-order delay differential equation that is representative of feedback systems with relays (switches) that actuate after a fixed delay. The system under study exhibits strong multirhythmicity, the coexistence of many stable periodic solutions for the same values of the parameters. We present a detailed study of these periodic solutions and their bifurcations. 
Starting from an integro-differential model, we show how to reduce the system to a set of finite-dimensional maps. We then demonstrate that the parameter regions of existence of periodic solutions can be understood in terms of discontinuity induced bifurcations and their stability is determined by smooth bifurcations. Using this technique we are able to show that slowly oscillating solutions are always stable if they exist. We also demonstrate the coexistence of  stable periodic solutions with quasiperiodic solutions.

 \end{abstract}

%\pacs{42.65.Sf,  05.45.-a }

% 05.45.Gg Control of chaos, applications of chaos
% 05.45.-a Nonlinear dynamics and chaos
% 05.45.Jn High-dimensional chaos
% 42.65.Sf Dynamics of nonlinear optical systems; optical instabilities; optical chaos and complexity; and optical spatio-temporal dynamics
% 89.75.Da Systems obeying scaling laws 
%\keywords{}

\maketitle

\section{Introduction}

This paper is concerned with the dynamics of time-delayed relay systems. 
Such systems have great practical importance, as relay control-systems are applied in many different areas of engineering.
In relay control, the control signal is a piecewise constant function of the measured output, typically switching between just two values. In addition, the control signal is often delayed due to finite signal transmission and processing times, sampling delays, or other latencies in control loops. As a result, the controlled system is described by a nonsmooth delay differential equation (DDE). 
Aside from control applications, nonsmooth DDEs arise as descriptions of naturally occurring systems for which nonlinearities are well approximated by functions that take on discrete values~\cite{Barton2006, Ryan2020}. 
Due to their importance, nonsmooth DDEs have been the focus of much recent attention~\cite{Fridman2000, Fridman2002, Shustin2003, Benadero2019, Sieber2006, Barton2006, Ryan2020}.

Another reason to study relay systems is that analytic results become possible.
Delay differential equations have, generically, an infinite dimensional state space \cite{Hale2002}, corresponding to the fact that  the initial condition consists of the history of the system over the entire delay interval.  The high-dimensionality makes the study of DDEs challenging.  Under relay feedback,  the dynamics can be reduced to the dynamics of finite dimensional maps, which leads to significant simplifications of the analytic treatment. The trade-off is that the interplay between the discontinuous relay feedback (switching) and delay leads to new-types of bifurcation scenarios, so-called discontinuity-induced bifurcations~\cite{Sieber2006}.

The dynamics of first order time-delayed relay systems is well understood. Under some mild assumption on the DDE, it can be shown that solutions will converge to an orbitally stable slowly oscillating periodic solution~\cite{Fridman2000,Fridman2002}. 

However, first order systems are insufficient to describe many naturally occuring phenomena and technological applications. Often second order models are required to capture the essence of the dynamics. Yet, the classification of possible solutions of second order time-delayed relay systems remains largely incomplete. Even information concerning the properties of just the periodic solutions is only partially answered, especially for arbitrary choices of parameters. 

For second order models, an interesting new possibility arises in control applications because the relay signal can now depend on two variables, e.g. the system position and velocity~\cite{Sieber2006}.

In this papers we study one of the simplest generic models, a linear second order time-delayed relay systems in which velocity serves as the feedback signal.
We show that the system demonstrates a complex bifurcation structure with significant multirhythmicity, the coexistence of multiple stable periodic solutions~\cite{Weicker2015}. Due to the combination of linear flow and piecewise-constant feedback, the system can be reduced to a set of finite dimensional maps. We use these maps to analytically determine regions of existence and stability of periodic solutions.

\section{Background}

Second order linear DDEs with time delayed relay feedback have of the form
\begin{align} \label{2ndorderDDE}
\alpha \, \ddot{y} +  \beta\, \dot{y}+ \gamma \, y = \sigma \, \mathrm{sign}[g(y(\mathsf{t}-\tau), \dot{y}(\mathsf{t}-\tau))],
\end{align}
where $\alpha, \beta, \gamma, \tau$ are real constants, the dot denotes the derivative with respect to time $\mathsf{t}$, and $\sigma= \pm 1$ signifies positive or negative feedback.
The function $g(y,\dot{y}):\mathbb{R}^2 \to \mathbb{R}$ divides the $y,\dot{y}$-plane into two domains, $\{g < 0\}$ and $\{g \ge 0\}$ and the relay-nonlinearity is modeled as 
\begin{align} \label{relaydef}
 \mathrm{sign}(x) = \begin{cases} +1 & x>0 \\ -1 & x < 0 \\ 0 & x=0 \end{cases}.
\end{align}

If parameters $\alpha, \beta$ are positive and $\gamma$ is negative, then \eq{2ndorderDDE} is analogous to an inverted pendulum with relay feedback.  The case where the delayed relay signal depends on the position, $g = y(\mathsf{t}-\tau)$, has been discussed in~\cite{Shustin2003}. It was shown that a unique stable periodic orbit is possible under certain conditions on the parameters. 
These conditions can be relaxed if the delayed relay signal depends not only position but velocity as well, as shown in \cite{Sieber2006}.

If parameters $\alpha, \beta, \gamma$ are positive, then \eq{2ndorderDDE} represents a harmonic oscillator with relay feedback.  
Such a model with $\beta=0$ and the delayed relay signal depending on the position is discussed in \cite{Heiden1990, Barton2006}, with the model originating from experimental studies on the pupil light reflex~\cite{Milton1990}. A rich set of bifurcations and solutions was identified~\cite{Barton2006}.

In this paper we discuss \eq{2ndorderDDE} with positive parameters $\alpha, \beta, \gamma$ and the delayed relay signal depending on the velocity. This situation was also considered in~\cite{Benadero2019} in an investigation of switch-mode power converters. Using numerical case studies, the authors find periodic orbits as well as smooth and discontinuity induced bifurcations similar to our results. In contradistinction to their work, we present explicit analytic solutions, which allow us to develop a more complete picture of this system.

\section{Model}

We consider a delayed feedback system consisting of a two-pole bandpass filter with scalar (dimensionless) output $x$ that is delayed by $\tau$ and passed through a relay, with the resulting signal serving as the scalar input to the bandpass filter. 

A convenient form for a two-pole bandpass-filter transfer-function, written in terms of angular frequency $\omega$, is
\begin{align}
H(\omega)  &= \frac{1}{1 + iQ\left({\omega}/{\omega_\text{c} }- {\omega_\text{c}}/{\omega}\right)} ,
\end{align}
where $\omega_\text{c}$ is the center angular frequency, i.e. the angular frequency of maximum transmission with $|H(\omega_\text{c})| = 1$, and $Q$ is the quality factor, defined as the ratio of the center frequency and filter bandwidth~\cite{Dimopoulos2012}.
The zero-response (or natural-response) of the filter describes the transient response of the filter due to nonzero initial conditions. It is oscillatory if $Q>1/2$, the underdamped regime, and nonoscillatory if $Q<1/2$, the overdamped regime.

We can associate the following integro-differential equation to the system
\begin{align} \label{integroDif1}
x + \frac{Q}{\omega_\text{c}} \dv{x}{\mathsf{t}} + Q \, \omega_\text{c} \int^{\mathsf{t}} \!\! x(s) \, \dd s =  \sigma \, \mathrm{sign}(x(\mathsf{t} - \tau)),
\end{align}
where the RHS is the input signal to the bandpass filter (LHS).
We distinguish positive ($\sigma=+1$) and negative ($\sigma=-1$) feedback.

For our analysis it is convenient to reference time to the delay by introducing dimensionless time $t$ via $t =  \mathsf{t} / \tau$ and to introduce the parameter $\Omega$, defined as the product of the filter's center frequency and the delay, 
\begin{align} \label{OmegaDef}
\Omega = \omega_\text{c} \, \tau.
\end{align} 
Furthermore, we introduce a variable $y$ that satisfies $\dot{y} = Q \Omega \, x$, where the dot denotes the derivative with respect to dimensionless time $t$.  
This yields the nonsmooth delay differential equation
\begin{subequations}  \label{DDEmodel}
\begin{align}
{Q}\,{\Omega}^{-1} \,  \dot{x} &= -  x - y + \sigma \,  \mathrm{sign}(x(t-1)) \\
\dot{y} &= Q \, \Omega \, x ,
\end{align}
\end{subequations}
which is equivalent to \eq{2ndorderDDE} with $g = \dot{y}(\mathsf{t}-\tau)$, $\alpha=(\tau/\Omega)^{2}$, $\beta=\tau/(Q \Omega)$, and $\gamma=1$.
Model~(\ref{DDEmodel})  depends on two positive dimensionless parameters: (1) the quality factor $Q$, set by the fractional bandwidth of the filter, and (2) the parameter $\Omega$. For a fixed center frequency, $\Omega$  increases proportional to the delay. Alternatively, for a fixed delay,  $\Omega$  increases with the filter's center frequency. If $\Omega = 2 \pi$, then a sinusoidal signal with a period equal to the delay $\tau$ has a frequency that coincides with the filter's center frequency. If  $\Omega > 2 \pi$ ($\Omega < 2 \pi$), then a sinusoidal signal with period $\tau$ has a frequency smaller (larger)  than the center frequency and will be attenuated by the filter.

\section{Sample Solutions and Symbolic Representation}

In order to solve the delayed relay system (\ref{DDEmodel}), we note that it is sufficient to keep track of the headpoint coordinates $x(t),y(t)$ and the sign of $x$ in the delay interval or, equivalently, the times $\tau_n$ at which $x$ crosses zero with $t-1< \tau_{k}< \tau_{k-1} < \ldots < \tau_{1} \le t$. 
The state of the system at time $t$ is, therefore, captured by a tuple of finite but variable length $k+2$, $(x,y;\tau_1,\tau_2,\ldots,\tau_k)$.
The solution can be obtained in terms of a discrete time map acting on the state that maps between key events.
There are two key events that occur:
\begin{enumerate}
\item  A \emph{zero} element is added at time $t=t_n$. When $x$ passes through zero, a zero element equal to $t_n$ is added to the state; increasing the length of the tuple by one. This event does not immediately cause the relay feedback to change but it will switch the sign of the feedback at time $t_n + 1$. Such an event is denoted by the symbol $Z$ if $x$ transitions from $x < 0$ to $x > 0$ and by $\overline{Z}$ for the opposite transition. 
\item A zero-crossing time is removed from the \emph{history}. When $\tau_{k}(t)=t-1$, $\tau_{k}$ is deleted from the state tuple.  The sign of the relay feedback switches. Such an event is denoted by the symbol $H$ if the feedback switch is due to a transition from $x(t-1) < 0$ to $x(t-1) > 0$ and by $\overline{H}$ for the opposite transition.
\end{enumerate}

Between consecutive events, the feedback term is constant, either $+1$ or $-1$. The evolution is given by the solution of a linear ordinary differential equation. Explicit solution of the ODE allows us to construct an iterative map that moves the system forward in time from one event to the next. 

\begin{figure}
\begin{center}
\includegraphics[width= .99 \columnwidth]{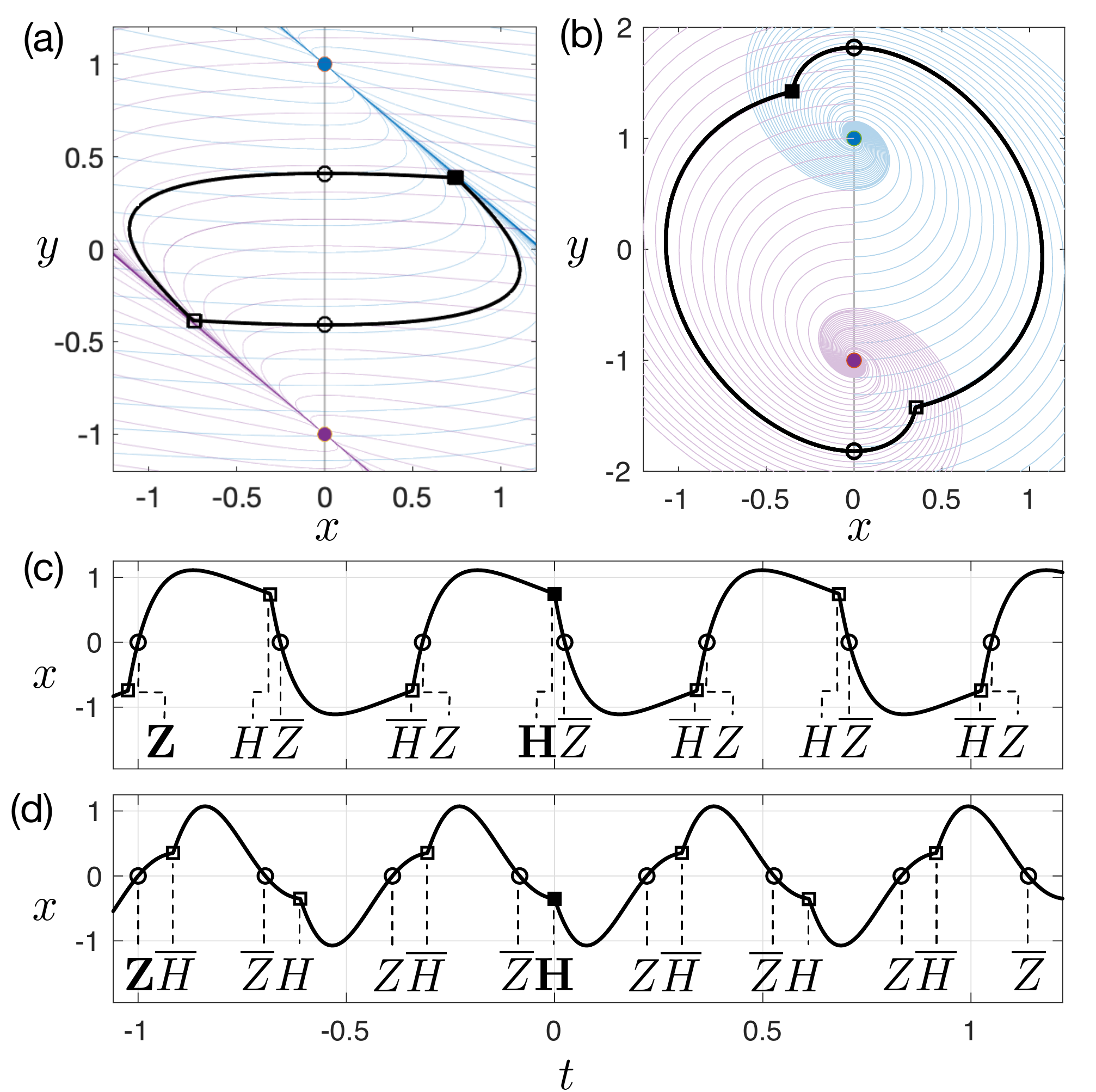}
\caption{Periodic solutions for negative feedback ($\sigma=-1$): (a),(b) projection onto  $x-y$ plane (black), ODE flows associated with a fixed constant sign of the relay term (blue and purple), and switching manifold (grey); (c),(d) $x(t)$ with events indicated. (a),(c) Parameters:  $Q=0.4$, $\Omega = 7$. Symmetric frequency-two solution $[H, \overline{Z}, \overline{H}, Z]_2^S$. (b),(d) Parameters: $Q=1.5$, $\Omega = 14$, symmetric frequency-three solution $[H, Z, \overline{H}, \overline{Z}]^S_3$.}
\label{fig:FlowPlot}
\end{center}
\end{figure}

The ODE flow for constant feedback has a stable fixed point that lies on the switching manifold and is either a node (see Fig.~\ref{fig:FlowPlot}a) or a spiral (see Fig.~\ref{fig:FlowPlot}b). In the latter case, trajectories are guaranteed to cross the switching manifold such that solutions of \dde{DDEmodel} are necessarily oscillatory. In the former case, there are initial conditions that result in non-oscillatory solutions to \dde{DDEmodel}, ones that approach the fixed point without ever crossing the switching manifold. In this paper we focus on oscillatory solutions and will exclude such initial conditions from consideration. We also exclude from consideration the unstable trivial solution ($x=y=0$).

To any oscillatory solution corresponds a symbolic-sequence representing the events. Periodic solutions are represented by a repeating sequence of events $[S_1, S_2, \ldots, S_n]$ with $S_i \in \left\{ H, \overline{H}, Z, \overline{Z}  \right\}$. Since every cyclic permutation of the repeating sequence represents the same solution, we start all sequences with $H$, for the sake of consistency.

To label distinct periodic solutions that have an identical symbol sequence, we define a discrete number -- the \emph{oscillation frequency} $\nu$, which is the number of $x$-zero crossings on the unit time-interval of the delay preceding a time $t$ at which $x(t)$ is zero. A periodic solution is said to be \emph{slowly oscillating} if $\nu=0$ and \emph{rapidly oscillating} if $\nu>0$.

We further label solutions by their symmetry. We note that \dde{DDEmodel} remains invariant under the operation $(x,y) \to (-x,-y)$. Periodic solutions that posses this symmetry are called symmetric, otherwise they are asymmetric. Asymmetric periodic solutions come in symmetry-related pairs.

Thus, our scheme for labeling periodic solutions symbolically is
\begin{align}
[ S_1, S_2, \ldots, S_n]_\nu^{\frak{s}}
\end{align}
with event symbols $S_i \in \left\{ H, \overline{H}, Z, \overline{Z}  \right\}$, frequency  $\nu = 0,1,2, \ldots$, and symmetry label $\frak{s} \in \left\{ S, A  \right\}$ for symmetric and asymmetric solutions, respectively.

As an example, we depict in Fig.~\ref{fig:FlowPlot} two periodic solutions and, for each, indicate $Z$, $\overline{Z}$ events by circles and $H$, $\overline{H}$ events by squares. The solution shown in Fig.~\ref{fig:FlowPlot}(c) is a symmetric solution with frequency $\nu =2$ and repeating four-symbol sequence $[H, \overline{Z}, \overline{H}, Z]_2^S$. The periodic solution in Fig.~\ref{fig:FlowPlot}(d)  is a symmetric $\nu =3$ solution with symbol-sequence $[H, Z, \overline{H}, \overline{Z}]_3^S$.

The frequency label $\nu$ can be related to the symbolic representation by noting that every $H$ ($\overline{H}$) event is associated with a $Z$ ($\overline{Z}$) event one delay time in the past. The frequency $\nu$ can be understood as the number of  $Z$/$\overline{Z}$ symbols in between. As an example, consider the $\nu=2$ solution in Fig.~\ref{fig:FlowPlot}(c). The $H$ event at $t=0$ that is indicated by the filled square is associated with the $Z$ event at time $t=-1$ and there is one $\overline{Z}$ and one $Z$ in between.

\section{Periodic Solutions}

A periodic solution of \eq{integroDif1} with period $\overline{P}$ and discrete frequency $\nu$ will also be a solution of \eq{integroDif1} if the delay is changed to $\tau' = n \, \overline{P} + \tau$ ($n=1, 2, \ldots$) because this mapping leaves \eq{integroDif1} unchanged. With respect to the delay $\tau'$, this periodic solution has a discrete frequency $\nu' > \nu$  (for a 4-symbol symmetric solution $\nu' = \nu + 2 n$). 
In terms of the dimensionless \dde{DDEmodel}, this mapping implies $x_{\nu'}(t | \Omega', Q)=x_\nu(t | \Omega, Q)$. Here, $x_\nu(t | \Omega, Q)$ denotes a periodic solution of \dde{DDEmodel} with frequency $\nu$ and  $x_{\nu'}(t | \Omega', Q)$ the corresponding $\nu'$-frequency  solution of \dde{DDEmodel} with the $\Omega$ parameter changed to $\Omega' = \Omega + \omega_\text{c} n \overline{P}$. 
Under this mapping, the $\Omega$-interval of existence of $x_\nu(t | \Omega, Q)$ results in a corresponding interval of existence of the $\nu'$-frequency solution. If these intervals overlap, the $\nu$ and $\nu'$ frequency solutions coexist, suggesting that \dde{DDEmodel} may have an infinite number of coexisting periodic solutions (since $n$ is arbitrary).
While suggestive, it is necessary to study solutions and their domain of existence in more detail in order to confirm coexistence and determine stability.
%If one considers now smooth changes of the parameter $\Omega'$ back to $\Omega$ and assumes that the $x_{\nu'}$ solution deforms but continues to exist as a $\nu'$-frequency solution, then \dde{DDEmodel} must have periodic solutions of frequencies $\nu$ and $\nu'$. The existence of a single periodic solution would then imply an infinite number of coexisting periodic solutions (since $n$ is arbitrary). Of course, the assumption in this argument can break down in many ways, such as for solutions that are destroyed by bifurcations upon variation of $\Omega$. It is therefore necessary to study solutions and their domain of existence in more detail. 

\begin{figure}
\begin{center}
\includegraphics[width= .99 \columnwidth]{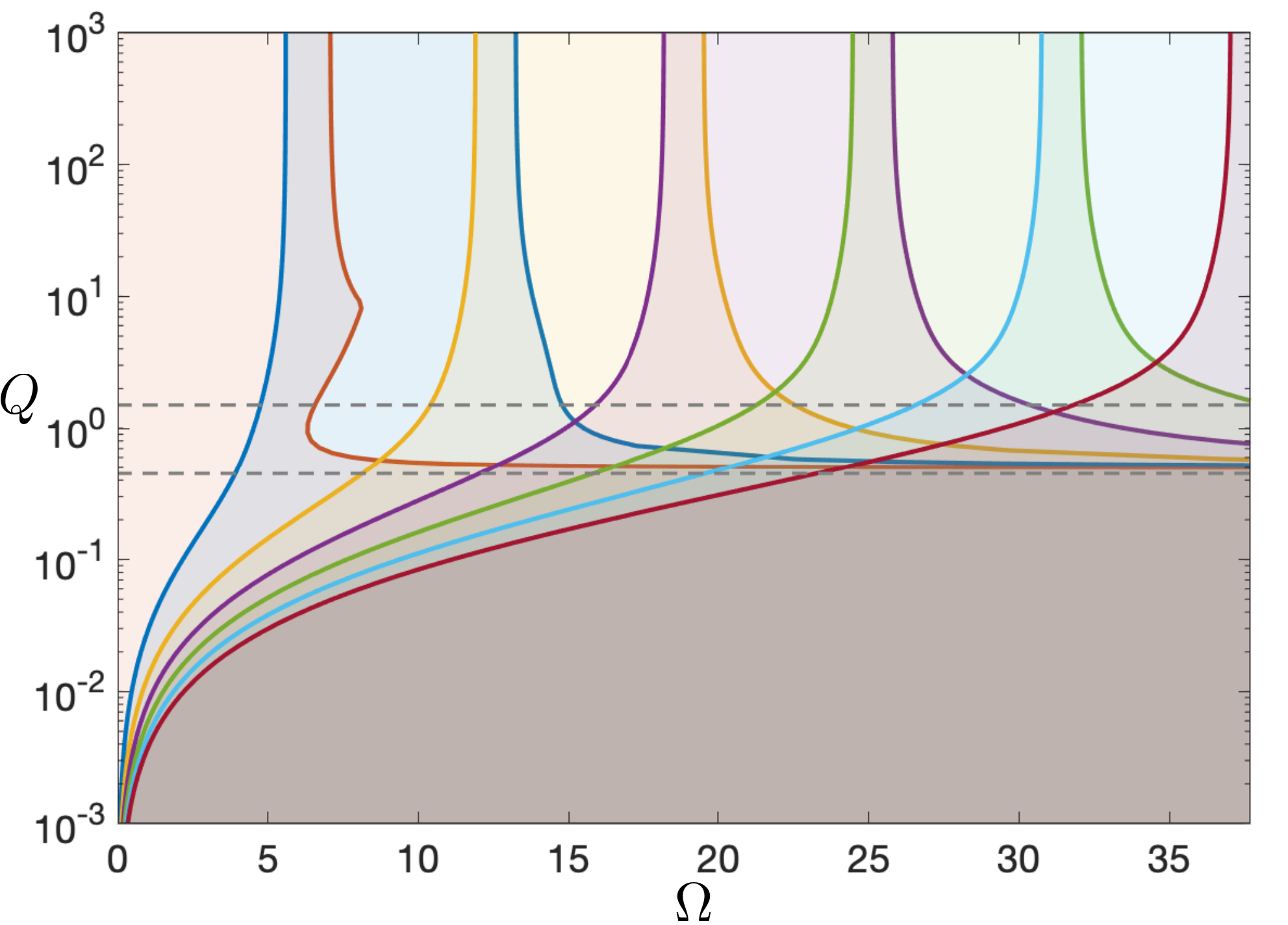}
\caption{Stable symmetric periodic solutions in $\Omega$--$Q$ parameter plane for negative feedback. Each colored region corresponds to a stable ``mode.'' At the boundaries (colored lines) the modes are destabilized through bifurcations. }
\label{fig:RegionPlot}
\end{center}
\end{figure}

\begin{figure}
\begin{center}
\includegraphics[width= .99 \columnwidth]{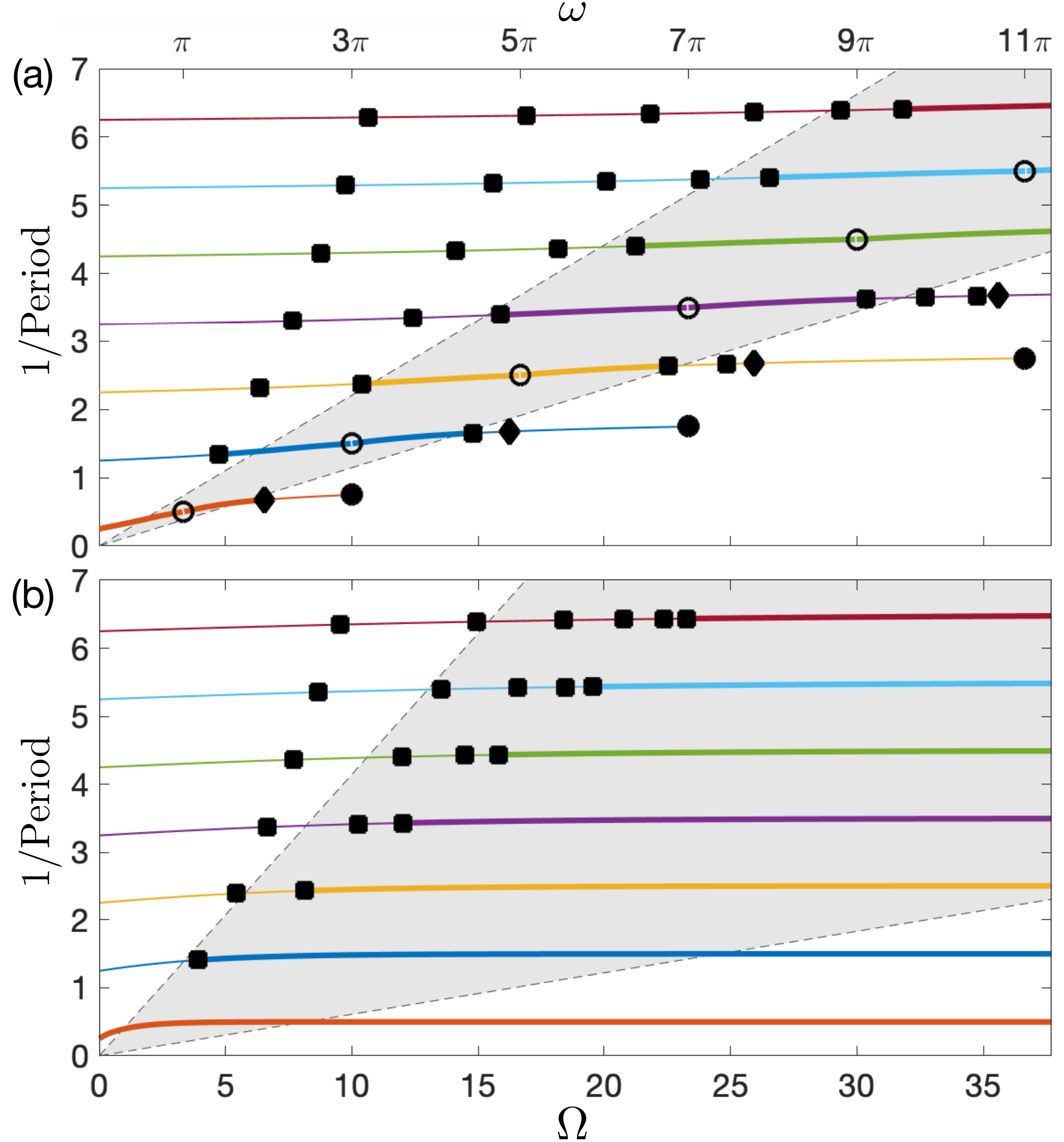}
\caption{Inverse period versus $\Omega$ for negative feedback and (a) underdamped regime with $Q=1.5$ (b) overdamped regime with $Q=0.45$ (shown by dashed lines in Fig.~\ref{fig:RegionPlot}). Periodic solutions: (thick line) stable,  (thin line) unstable. Smooth bifurcations:  (squares) Neimark-Sacker bifurcation, (diamonds) pitchfork bifurcation. Discontinuity induced transitions: (circles) the $H/\overline{H}$ headpoint is at switching manifold, i.e. $x_H = 0$. Filter: (shaded region) passband as defined by 3dB frequencies. 
}
\label{fig:PeriodPlotNF}
\end{center}
\end{figure}

Figures~\ref{fig:RegionPlot} and \ref{fig:PeriodPlotNF} demonstrate not only the presence of coexisting periodic solutions, but strong multirhythmicity of the system. 

In Fig.~\ref{fig:RegionPlot} each color represents a mode of stable periodic oscillations, where by a mode we mean a periodic solution of \dde{DDEmodel} that varies smoothly with parameters $Q,\Omega$. 
As seen by the overlapping regions in Fig.~\ref{fig:RegionPlot}, there are many stable periodic solutions that coexist. 
 We note that only the seven modes with the lowest frequency are shown and additional overlaps would appear as more modes are included.

In the underdamped regime ($Q>1/2$) holding $Q$ fixed and changing $\Omega$, periodic solutions are stable for a finite $\Omega$ interval. This is consistent with the intuition that the frequency of each oscillating mode should be commensurate with the filter's passband. Surprisingly, this intuition fails for the overdamped regime ($Q<1/2$). 

In the overdamped regime, the number of coexisting stable modes increases without bound as $\Omega$ increases, as is indicated in Fig.~\ref{fig:RegionPlot}, where the overlap of all colors (resulting in brown) means that all the depicted modes are stable. For this reason,  multirhythmicity is particularly pronounced.
 
To demonstrate the strong  multirhythmicity in the overdamped regime more clearly, we show in Fig.~\ref{fig:PeriodPlotNF}b the inverse period as a function of $\Omega$ for the seven periodic solutions $[H, \overline{Z}, \overline{H}, Z]_{\nu}^S$ with $\nu = 0, 2, 4, 6, 8, 10, 12$ (bottom to top) that correspond to the seven modes shown in Fig.~\ref{fig:RegionPlot}.
As $\Omega$ increases, all periodic solutions are stabilized via torus bifurcations (Neimark-Sacker bifurcations of the Poincar\'{e} map)  and remain stable. 
As seen in Fig.~\ref{fig:PeriodPlotNF}b, the periodic solutions remain stable even if their fundamental frequency is well below the low-frequency 3dB cut-off  of the filter. 
We argue that no additional bifurcations occur for values of $\Omega$ beyond those shown in Fig.~\ref{fig:PeriodPlotNF}b (see Sec.~\ref{sec:theory}). This, for a fixed filter and upon recalling the definition of $\Omega$ given by \eq{OmegaDef}, implies that the number of coexisting stable modes continues to grow as the delay $\tau$ is increased.

In contradistinction to the overdamped regime, in the underdamped regime the periodic oscillations gain and loose their stability as $\Omega$ increases. As seen in Fig.~\ref{fig:PeriodPlotNF}a, there is a strong correspondence between the intervals of $\Omega$ in which modes are stable and the passband of the filter.
The depicted symmetric periodic solutions gain and loose their stability via Neimark-Sacker bifurcations, the only exception is the lowest frequency mode, shown as the red bottom most line in  Fig.~\ref{fig:PeriodPlotNF}a, which loses stability in a pitchfork bifurcation.

In Fig.~\ref{fig:PeriodPlotNF}a, circles indicate for each mode the values of $\Omega$ at which the $H/\overline{H}$ headpoint is at the switching manifold, i.e. $x_H=0$.  At the first such point (open circles) the headpoint simply passes through the switching manifold and the periodic solution continues to exist.  However, the symbol sequence switches: 
$[H, \overline{Z}, \overline{H}, Z]_{\nu}^S$ with $\nu = 0, 2, 4, 6, 8, 10$ (bottom to top) becomes $[H, Z, \overline{H}, \overline{Z}]_\nu^S$ with $\nu = 1, 3, 5, 7, 9, 11$, respectively.
For example, the curve second from the bottom in Fig.~\ref{fig:PeriodPlotNF}a (blue curve) is a single mode associated with periodic solution labels $[H, \overline{Z}, \overline{H}, Z]_{2}^S$ and $[H, Z, \overline{H}, \overline{Z}]_3^S$. We refer to this mode as the $\nu=2,3$ four-symbol symmetric mode.
At the second $x_H=0$ point of each mode (filled circles), there is a discontinuity induced bifurcation and the mode ceases to exist.

\section{Bifurcation Analysis \label{sec:numerics}}

Next we turn to a detailed numerical bifurcation analysis of one of the modes, the $\nu=2,3$-mode, in order to elucidate 
 stability regions, bifurcations, and boundaries of existence. The bifurcation diagram produced is representative, we find similar bifurcations for other modes. 

\begin{figure}
\begin{center}
\includegraphics[width= .99 \columnwidth]{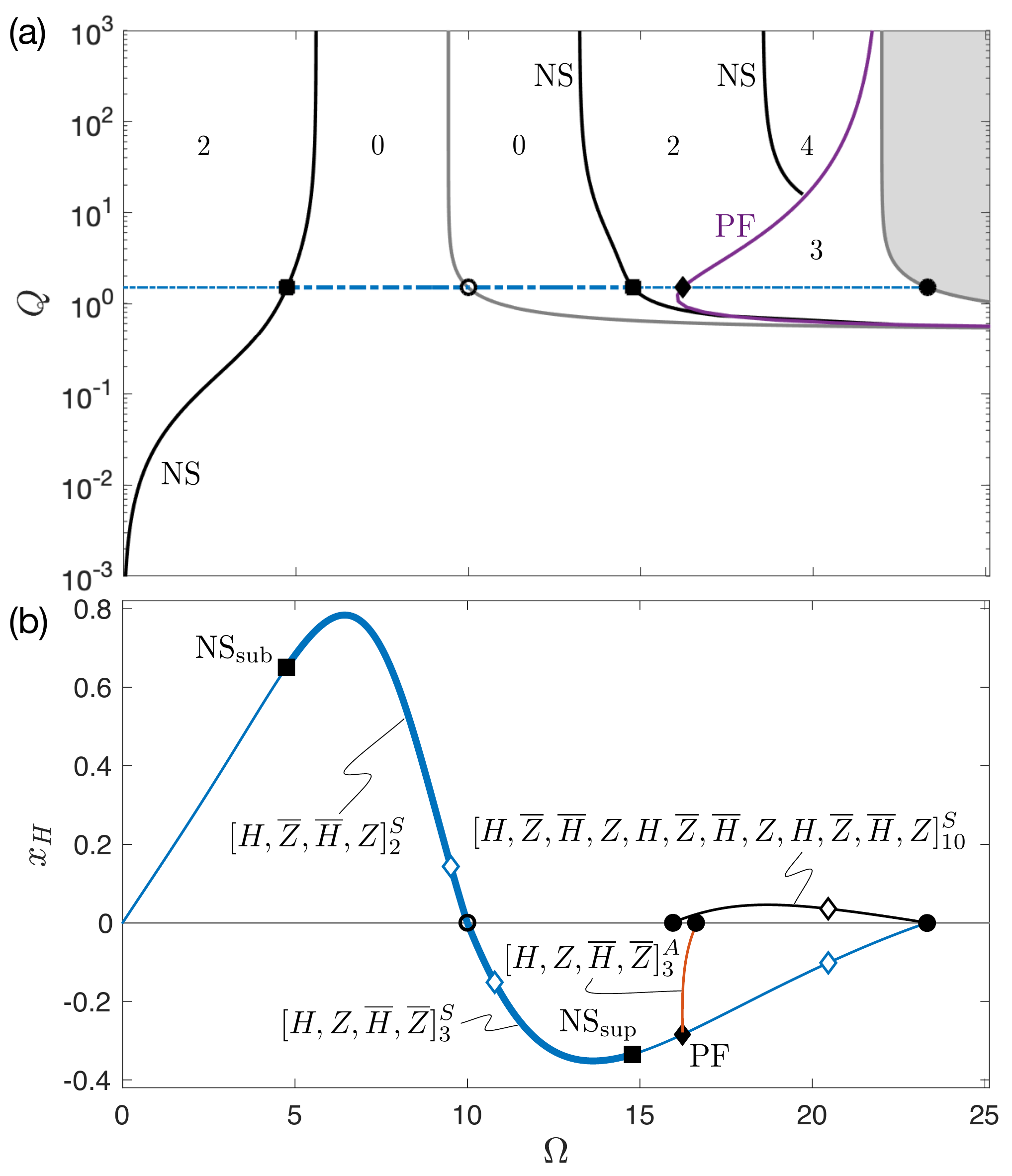}
\caption{Bifurcations of the $\nu=2,3$-mode for negative feedback. (a) $\Omega$--$Q$ parameter plane (b) $x_H$ of $\nu=2,3$-mode  (blue) versus $\Omega$ for $Q=1.5$ with connecting asymmetric (red) and 12-symbol (black)  branches also shown. (Thick lines) stable, (thin lines) unstable.
Smooth bifurcations: Neimark-Sacker (filled squares) and  pitchfork (filled diamond). Discontinuity induced transitions/corner collisions (circles). 
(Open diamonds) see Fig.~\ref{fig:BifPlotk2k3_Extc}.
}
\label{fig:BifPlotk2k3}
\end{center}
\end{figure}

For each symbol sequence and given discrete frequency $\nu$, the DDE reduces to a map with fixed dimension. The fixed points of this map correspond to periodic orbits. Orbit location and stability can be determined numerically by using the MatContM continuation software~\cite{Govaerts2008} as well as analytically (see Sec.~\ref{sec:theory} for details).

In Fig.~\ref{fig:BifPlotk2k3}, we label bifurcations by the fixed point bifurcations of the map. The Neimark-Sacker (NS) and pitchfork bifurcation (PF) curves in Fig.~\ref{fig:BifPlotk2k3}a are analytic and coincide with curves obtained using numeric continuation. Also given is the number of unstable directions, i.e. the number of eigenvalues of the Jacobian of the map with magnitude larger one.
The parameter region with zero unstable directions, including the two bounding Neimark-Sacker bifurcation curves, corresponds to the region enclosed by blue curves in Fig.~\ref{fig:RegionPlot}.

Fixing the parameter $Q$ to $Q=1.5$, we obtain the curves in Fig.~\ref{fig:BifPlotk2k3}b via numeric continuation. Shown is the value of $x$ at the $H$ event as a function of $\Omega$. It is seen that the $\nu=2,3$-mode (blue) is stable (thick line) in between the subcritical Neimark-Sacker bifurcation at $\Omega = 4.75$ and the supercritical  Neimark-Sacker bifurcation at $\Omega = 14.78$. All other periodic solutions shown are unstable. 

The conditions for construction of the finite dimensional maps that capture solutions of the DDE are that the projection of the solution onto the $x$-$y$-plane has the following two properties:
\begin{enumerate}
\item  The set of $Z/\overline{Z}$-event points is finite and disjoined from the set of $H/\overline{H}$-event points.
\item The flows are transverse to the switching manifold at all intersections.
\end{enumerate}
If one of those conditions is violated a discontinuity induced transitions occurs. 

Due to the symmetry of \dde{DDEmodel},  the flow is always transversal to the switching manifold, with the consequence that all  discontinuity induced transitions arise due to a violation of condition (1).

\begin{figure}
\begin{center}
\includegraphics[width= .99 \columnwidth]{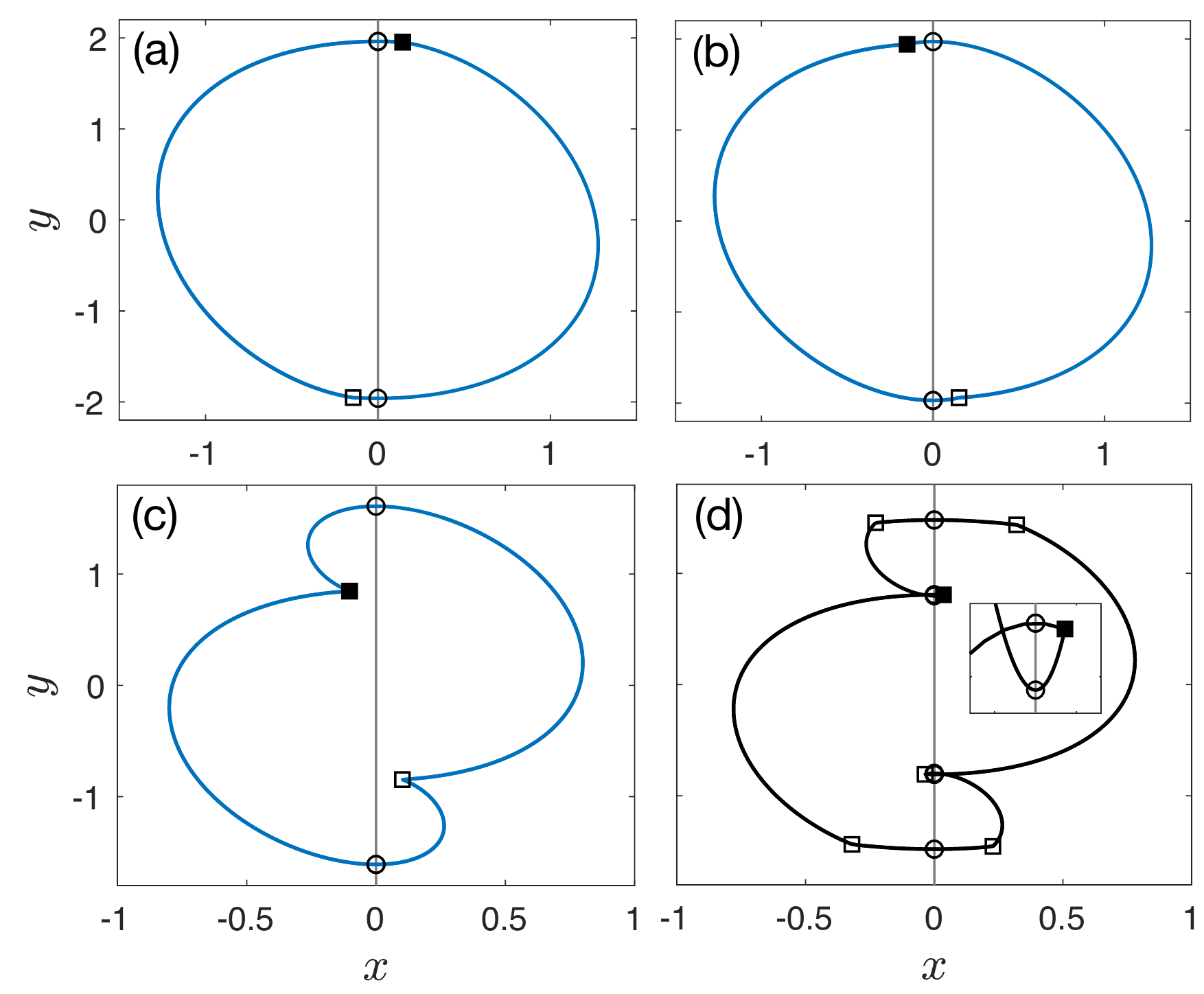}
\caption{  Solutions corresponding to open diamonds in Fig.~\ref{fig:BifPlotk2k3}b with $H$-event used for $x_H$ indicated by a solid square. 4-symbol symmetric $\nu=2,3$-mode: (a) $\Omega = 9.51$, $\nu=2$  (b) $\Omega = 10.79$, $\nu=3$  (c) $\Omega = 20.46$, $\nu=3$.  12-symbol symmetric  solution:  (d)  $\Omega = 20.46$, $\nu=10$.
Inset: detail showing that $x_H$ is positive (solid square) and that there are additional $x$-zero crossings close by (open circles).
}
\label{fig:BifPlotk2k3_Extc}
\end{center}
\end{figure}

The violation of condition (1) means that an $H/\overline{H}$-point of a solution reaches the switching manifold. Since at an $H/\overline{H}$-point the solution switches from one to the other flow, the graph of the solution typically has a corner. For this reason, such discontinuity induced transitions are called \emph{corner collisions}. We should mention that the corner collisions in our system violate one of the genericity conditions of \citet{Sieber2006} (condition 10b in \cite{Sieber2006}) due to the symmetry and perfect linearity of our system. Nevertheless, the corner collisions in our system do fall into the typical two main types:

In the first type, the part of the periodic orbit that is close to the colliding $H/\overline{H}$-point intersects the switching manifold transversally.  That is, in the vicinity of the $H/\overline{H}$-point both flows cross the switching manifold in the same direction.  In this case, as the bifurcation parameter is changed, the headpoint moves smoothly through the switching manifold, resulting in a smooth deformation of the periodic solutions. The mode, understood as a periodic solution of the DDE, continues to exist. This is shown in  Fig.~\ref{fig:BifPlotk2k3_Extc}a and Fig.~\ref{fig:BifPlotk2k3_Extc}b. However, because the symbol sequence switches and the number of zeros within the delay interval increases, the corresponding map changes discontinuously. Its dimension increases by one.

In the second type of corner collision, the part of the periodic orbit that is close to the colliding $H/\overline{H}$-point lies entirely on one side of the switching manifold for bifurcation parameter values slightly below the critical value. That is, in the vicinity of the $H/\overline{H}$-point the two flows cross the switching manifold in opposite directions.  In this case, the mode disappears in a discontinuity induced bifurcation.  This is shown in Fig.~\ref{fig:BifPlotk2k3_Extc}, where the 4-symbol symmetric $\nu=3$ mode is shown in Fig.~\ref{fig:BifPlotk2k3_Extc}c and for the same value of $\Omega$ a coexisting 12-symbol symmetric $\nu=10$ periodic solution is depicted in Fig.~\ref{fig:BifPlotk2k3_Extc}d. These two solutions coincide and cease to exist at the critical value of $\Omega$.

%
% FIGURE
%
\begin{figure}
\begin{center}
\includegraphics[width= .99 \columnwidth]{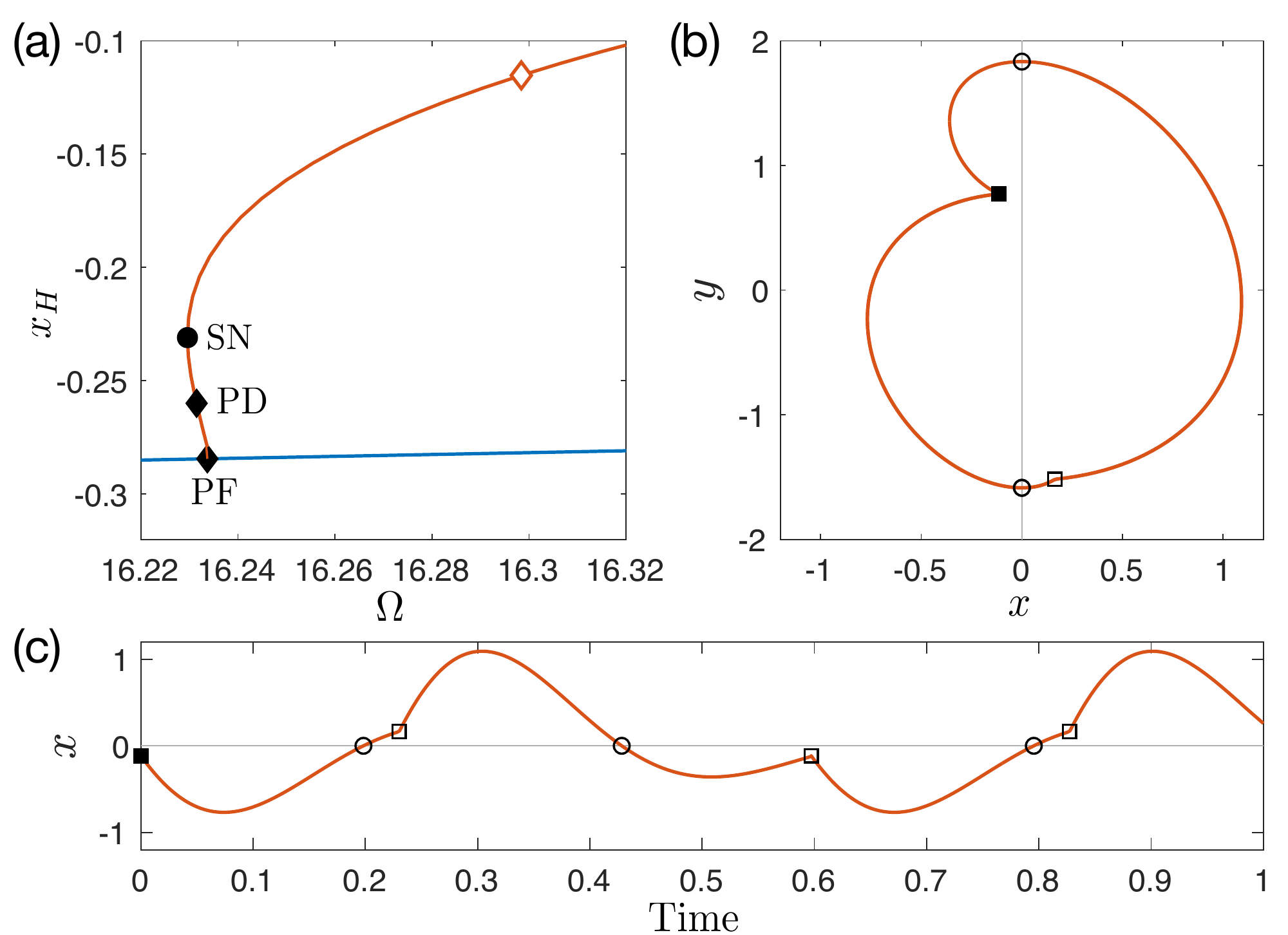}
\caption{Detail of asymmetric solution $[H,Z,\overline{H}, \overline{Z}]^A_3$ created via a pitchfork bifurcation (PF). (a) Zoom of Fig.~\ref{fig:BifPlotk2k3}b near PF: (PD) period doubling bifurcation, (SN) Saddle-Node bifurcation, (open diamond) solution shown in (b) and (c). 
}
\label{fig:BifPlotk2k3_Asym}
\end{center}
\end{figure}

Whereas  discontinuity induced transitions of periodic solutions of the DDE determine the range of validity of the corresponding map of a particular fixed dimension, standard smooth bifurcations of the periodic solutions of the DDE correspond to standard bifurcations of fixed points of the map. 

In Fig.~\ref{fig:BifPlotk2k3_Asym} we show that the pitchfork bifurcation of the $\nu=3$ symmetric solution leads to a symmetry related pair of asymmetric solutions (we only depict the branch of one of the two solutions). The asymmetric solution undergoes its own bifurcations, as seen in Fig.~\ref{fig:BifPlotk2k3_Asym}a. Although these bifurcations reduce the number of unstable eigenvalues, the solution remains unstable. The asymmetry is apparent in the $x-y$ projection of Fig.~\ref{fig:BifPlotk2k3_Asym}b and the timetrace depicted in Fig.~\ref{fig:BifPlotk2k3_Asym}c.
Similarly to the symmetric solutions, when the $H$-point of the asymmetric solution reaches the switching manifold ($x_H \to 0$), the solution disappears in a discontinuity induced bifurcation (see Fig.~\ref{fig:BifPlotk2k3}b).

We note that asymmetric solutions can be stable, such as the asymmetric solutions created via the $\nu=0,1$-mode's pitchfork-bifurcation that is shown as a diamond symbol on the lowest (red) curve in Fig.~\ref{fig:PeriodPlotNF}a.

\section{Quasiperiodic Solutions}

\begin{figure}
\begin{center}
\includegraphics[width= .99 \columnwidth]{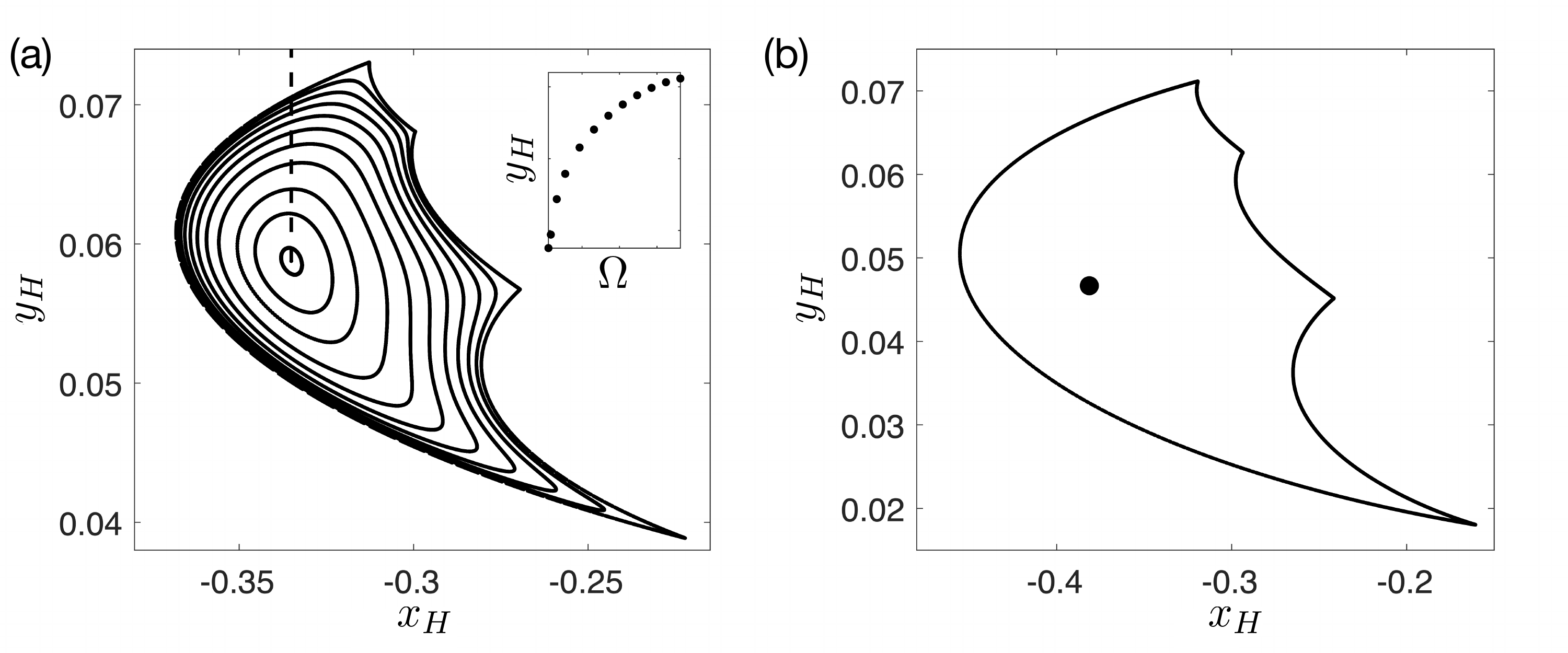}
\caption{ Quasiperiodic solutions for negative feedback ($\sigma=-1$). Projection onto  $x_H-y_H$ plane.  (a) Supercritical Neimark-Sacker bifurcation of 4-symbol symmetric $\nu=3$ periodic solution (filled square in Fig.~\ref{fig:BifPlotk2k3}). Parameters:  $Q=1.5$ and $\Omega$ ranging from $14.78$ to $14.84$.  Inset: $y_H$ along the dashed line as a function of $\Omega$ .  (b) Coexisting stable 4-symbol symmetric $\nu=3$ periodic solution and stable quasiperiodic solution for $Q=1.93$ and $\Omega = 14.56$. }
\label{fig:TorusPlot}
\end{center}
\end{figure}

Supercitical Neimark-Sacker bifurcations suggest the existence of stable quasiperiodic solutions of \dde{DDEmodel}. We indeed find such solutions, as shown in Fig.~\ref{fig:TorusPlot}a, where we plot $10^5$ iterates of the discrete map between $H$ events for $Q=1.5$ and several values of $\Omega$.  Initial iterates associated with the transient approach of the attractor were discarded. For values of the parameter slightly larger than the supercritical Neimark-Sacker bifurcation value ($\Omega =14.78$), the map iterates form a closed curve, indicating the existence of a torus attractor of \dde{DDEmodel}. 
The size of the torus grows smoothly as $\Omega$ is increased, as shown in the inset of Fig.~\ref{fig:TorusPlot}a. 
Numerically we were unable to find a stable torus attractor for values greater than $\Omega =14.84$. It would be interesting to determine the cause and explain the peculiar attractor shape at the largest value of $\Omega$, but we did not pursue this question. 

By using previous solutions as initial conditions and changing both $Q$ and $\Omega$ by small increments, it is possible to follow the largest torus attractor to $Q=1.93$ and $\Omega = 14.56$. For these parameters, the torus attractor coexists with the stable periodic solution, which we show in Fig.~\ref{fig:TorusPlot}b. For fixed $Q=1.93$, stable torus attractors of growing amplitude are still created via a supercritical Neimark-Sacker bifurcation of the periodic solution as $\Omega$ is increased past the bifurcation value. However, the torus curve then ``folds back,'' such that large amplitude torus attractors exist for $\Omega$ values smaller than the bifurcation value, leading to the coexistence depicted in Fig.~\ref{fig:TorusPlot}b.

Thus, we find rich dynamics  in the underdamped regime. Not only are multiple stable periodic solutions present (multirhythmicity) but  there are also parameter regions where, in addition, stable quasiperiodic solutions coexist. 

In the overdamped regime,  all of the Neimark-Sacker bifurcations found are subcritical. We were not able to locate any stable torus attractors but cannot exclude their existence. What we find is strong multirhythmicity, a large number of stable coexisting periodic solutions for large $\Omega$. 

%
% FIGURE
%

\section{Theory \label{sec:theory}}

The DDE (\ref{DDEmodel}) is reduced to discrete-time finite dimensional mappings between events by exploiting the linearity of the system to explicitly calculate the flow for times in between events. Defining the headpoint at time $t=0$ with $\mathbf{v} = (x,y)^T$, the flow is given by 
\begin{align}
\Phi_+(t, \mathbf{v}) = \mathbf{A}(t) \, \mathbf{v} + \mathbf{b}(t)
\end{align}
if $\sigma \cdot \mathrm{sign}(x(-1)) =+1$ and by
\begin{align}
\Phi_-(t, \mathbf{v}) = \mathbf{A}(t) \, \mathbf{v} - \mathbf{b}(t).
\end{align}
 if $\sigma \cdot \mathrm{sign}(x(-1)) =-1$.
Here,
\begin{align} \label{Adef}
\mathbf{A}(t) =
 e^{- \mu t}  \begin{pmatrix}
{\scriptstyle \cos(\omega t)  -} \frac{\mu \sin(\omega t)}{\omega} & -  \frac{2 \mu \sin(\omega t)}{\omega} \\
 \frac{\mu ^2+\omega ^2}{2 \mu} \frac{\sin(\omega t)}{\omega} &  {\scriptstyle \cos(\omega t)+} \frac{\mu \sin(\omega t)}{\omega} \\
 \end{pmatrix} 
\end{align}
and
\begin{align}\label{bdef}
\mathbf{b}(t) =
\begin{pmatrix}
 e^{-\mu t} \,  2 \mu \, \sin(\omega t)/\omega   \\
 1 - e^{-\mu t} \left( \cos (\omega t) + \mu \sin(\omega t)/\omega \right)
  \end{pmatrix} 
\end{align}
and we introduced the abbreviations
\begin{align}  \label{muomdef}
\mu &= \frac{\Omega}{2 Q}  
&
\omega &=  \Omega \, \frac{\sqrt{4 Q^2 -1}}{2 Q}. 
\end{align}
%
%\begin{subequations} \label{muomdef}
%\begin{align} 
%\mu &= \frac{\Omega}{2 Q}  \label{mudef} \\ 
%\omega &= \frac{\sqrt{4 Q^2 -1}}{2 Q} \, \Omega. \label{omdef}
%\end{align}
%\end{subequations}
The damping constant $\mu$ is positive definite, whereas $\omega$ is positive real if $Q>1/2$ but imaginary if $Q<1/2$. In the latter case,  $\omega = i |\omega| = \Omega {\sqrt{1 -4 Q^2}}/({2 Q})$ and in \eq{Adef} and \eq{bdef} one may use the identities $\cos (\omega t) = \cosh(|\omega| t)$ and $\sin(\omega t)/\omega = \sinh(|\omega| t)/|\omega|$.
The flow $\Phi_\pm$ has a single stable fixed point at $\mathbf{v}^* = (x,y)^T = (0, \pm 1)$, which is a spiral if $Q>1/2$ (underdamped regime) and a node if $Q<1/2$ (overdamped regime).

As an example of a mapping between events, let us consider a time $t_n$ at which there is a $Z$-type event, such that $x(t_n)=0$, and assume $\sigma x(t_n-1) < 0$ as well as $\nu>0$, meaning that there is at least one zero crossing in the history interval. Whether the next event is of $H$-type or $Z$-type is determined by evaluating the time required to reach either event under the assumption that the feedback sign does not switch and then picking the event that occurs first.
The time interval to the subsequent $H$-type event is  
\begin{align}
 \delta = \tau_\nu +1 - t_n,
\end{align}
whereas the time interval $z$ to the next $Z$-type event is determined by solving
\begin{align}
\begin{pmatrix} 0 \\ y \end{pmatrix}= \Phi_-(z, \mathbf{v}_Z) 
\end{align}
for the smallest positive $z$. Here, $\vb{v}_{Z}=(0, y_{Z})^T$ is the $Z$-headpoint at time $t_n$.  
The map to the headpoint of the next event is then 
\begin{align}
 \mathbf{v} = \Phi_-(\min\{\delta,z\}, \mathbf{v}_{Z}). 
\end{align}
In addition to updating the headpoint, the times of zero crossings in the history interval need to be updated. 
If an $H$-type event follows the $Z$-type event ($\delta<z$), then the number of zero crossings is unchanged because one zero crossing is removed and another added. If a $Z$-type event follows ($z<\delta$), then the dimension of the state-vector is increased by one.
Similarly, if two $H$-type events follow one another, then the number of zero crossings and the dimension of the state-vector are reduced by one. 

Stable periodic solutions of the DDE (\ref{DDEmodel}) that  were found numerically are all of the same type; they all are symmetric 4-symbol solutions consisting of alternating $H$-type and $Z$-type events. As stable solutions are important for applications, we provide next details about the relevant Poincar\'{e} map, its fixed points and their bifurcations.

\subsection{Map of four symbol symmetric solutions}

We consider 4-symbol periodic solutions of \dde{DDEmodel} with a symbol sequence consisting of alternating $H$-type and $Z$-type events, that is, either a repetition of the sequence $H, Z, \overline{H}, \overline{Z}$ or of $H, \overline{Z}, \overline{H}, Z$. These solutions correspond to fixed points of a  Poincar\'{e}  map $\mathcal{P}_\nu$ that maps the solution forward by 4-events per step. 

Since these solutions have alternating $H$-type and $Z$-type events, the discrete frequency $\nu$, which counts the number of zero crossings in the history, remains constant.  

Although we found it advantageous to keep track of $H$ head-points and their position relative to the switching manifold when doing numerics, in terms of  theory it is convenient to map between $Z$-type events, in which case the  $\nu+1$ independent variables can be chosen to be the $\nu$ time intervals between $x$-zero-crossings and the $y$ coordinate of the headpoint. The $x$ coordinate of the $Z$-event headpoint is zero by definition.

\begin{figure}
\begin{center}
\includegraphics[width= .9 \columnwidth]{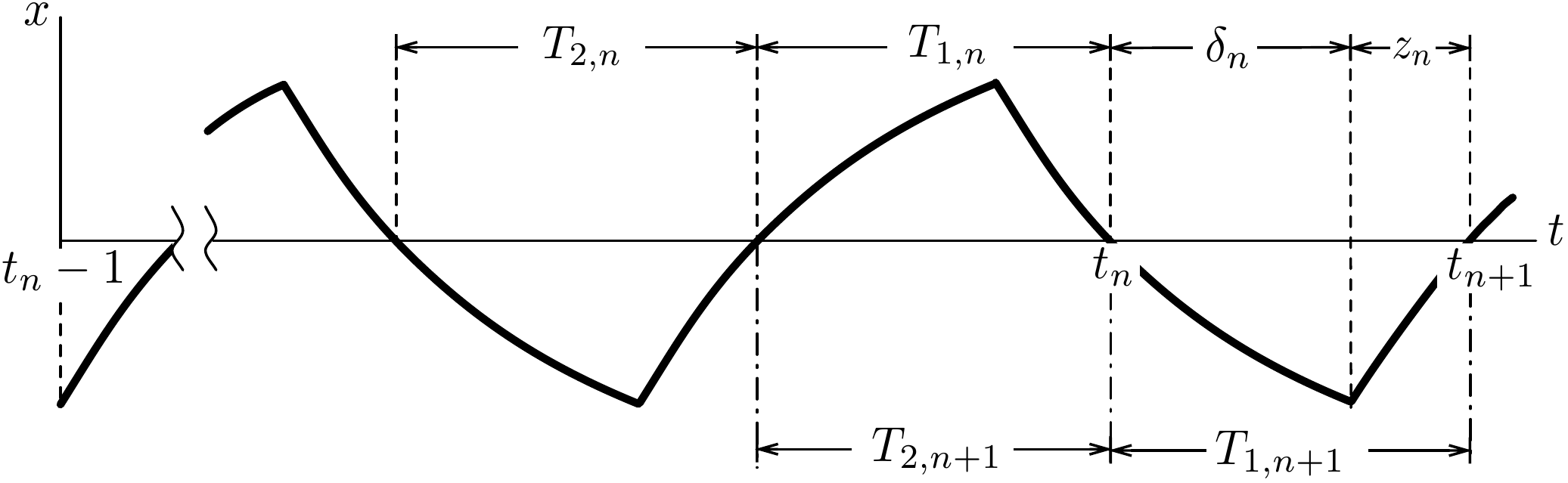}
\caption{Time intervals $T_{j,n}$ ($j=1 \ldots \nu$), $\delta_n$, and $z_n$.  Also shown, the time intervals after one iteration of the state vector, $\mathbf{s}_{n+1} = \mathcal{M}_\nu^+(\mathbf{s}_n)$.  }
\label{fig:MapSetup}
\end{center}
\end{figure}

In particular, let us consider a solution at some $Z$-type instance $t=t_n$, meaning $x(t_n)=0$. Furthermore,  assume that $\sigma x(t_n-1) < 0$ and there are $\nu$ zero crossings of $x$ of the history function at times $\tau_j$ with $j=1 \ldots \nu$ and $\tau_j \in (t_n-1, t_n)$. Let  $T_{j,n}$ denote the $\nu$ time intervals between $x$-zero-crossings, i.e. $T_{1,n} = t_n - \tau_1$ and $T_{j,n} = \tau_{j-1} - \tau_j$ for $j=2, \ldots, \nu$ (see Fig.~\ref{fig:MapSetup}).  Denote the headpoint time $t_n$ by  $\mathbf{v}_{Z,n} = (0 ,y_{Z,n})^T$ and define the $(\nu+1)$-dimensional state vector as
%
% FIGURE
%
\begin{align}
\mathbf{s}_n = (y_{Z,n}, T_{1,n},T_{2,n}, \ldots , T_{\nu,n})^T.
\end{align}
The headpoint of the subsequent $Z$-type event is given by the map
\begin{align} \label{headpointmap}
\mathbf{v}_{Z,n+1} = \Phi_+(z_n, \Phi_-(\delta_n, \mathbf{v}_{Z,n})).
\end{align}
The flow $\Phi_-$ shifts the solution until the next sign change of $x(t-1)$, which is an $H$ event for positive feedback ($\sigma = 1$) and an $\overline{H}$ event for negative feedback ($\sigma = -1$), with the time interval to the crossing being
\begin{align}  \label{deltandef}
\delta_n = 1 - \sum_{j=1}^\nu T_{j,n} = 1 - \sum_{i=2}^{\nu+1} s_{i,n},
\end{align}
if $\nu>0$ and $\delta_n=1$ if $\nu=0$.
The flow $\Phi_+$ then shifts the solution to the subsequent $Z$-type event. The time interval for this mapping is
\begin{align} \label{zndef}
z_{n}(\mathbf{s}_n) = \frac{1}{\omega} \arctan\left( \frac{ e^{-\mu \delta_{n} } \sin(\omega  \delta_{n}) (y_{Z,n} +1) }{ 2 - e^{-\mu \delta_{n} } \cos(\omega  \delta_{n}) (y_{Z,n} +1) }\right).
\end{align}

Utilizing the head-point mapping given by \eq{headpointmap}, we define a map $\mathcal{M}_\nu^+$ that updates the state vector, implementing the two-symbol $Z$-to-$Z$ shift, 
\begin{align}
\mathcal{M}_\nu^+(\mathbf{s}_n) = (y_{Z,n+1}, \delta_n+z_n, T_{1,n}, \ldots , T_{\nu-1,n})^T.
\end{align}
The map  $\mathcal{M}_\nu^+$ is \textit{frequency preserving} because one zero is removed and one zero is added to the history. In addition, denote with $\mathcal{M}_\nu^-$ the corresponding two-symbol shift  map which advances the solution through one more $H$-type event to the subsequent $Z$-type event.

The symmetry of \dde{DDEmodel} means that one can express $\mathcal{M}_\nu^-$ in terms of $\mathcal{M}^+_\nu$ by defining the operation of a sign flip as 
 \begin{align}
\mathcal{R}_\nu(\mathbf{s}_n) \equiv  \vb{R} \, \vb{s}_n 
\end{align} 
with $\vb{R}$ being a $(\nu+1)\times (\nu+1)$ diagonal matrix with $R_{11}=-1$ and $R_{jj}=1$ $(j = 2, \ldots, \nu+1)$. Then
\begin{align}
\mathcal{M}_\nu^- = \mathcal{R}_\nu \circ \mathcal{M}^+_\nu \circ \mathcal{R}_\nu.
\end{align} 
Therefore, the Pointcar\'{e} map $\mathcal{P}_\nu$ that maps a four-symbol $\nu$-frequency solution forward by four symbols is
\begin{align}
\mathcal{P}_\nu = \mathcal{R}_\nu \circ \mathcal{M}^+_\nu \circ \mathcal{R}_\nu \circ \mathcal{M}^+_\nu = (\mathcal{R}_\nu \circ \mathcal{M}^+_\nu)^2.
\end{align}
To investigate symmetric periodic solutions it suffices to study the fixed points of the map $\mathcal{M}_\nu = \mathcal{R}_\nu \circ \mathcal{M}^+_\nu$. Explicitely, $\mathbf{s}_{n+1} = \mathcal{M}_\nu(\mathbf{s}_n)$ is 
\begin{align} \label{FreqPresMapCompl}
\begin{split}
s_{1,n+1}&=  - 1 +2 \cos(\omega z_n)  \, e^{-\mu z_n} -   \\
& \quad -(s_{1,n} + 1) \, \cos(\omega [\delta_n + z_n]) \,e^{-\mu [\delta_n + z_n]}  \\
s_{2,n+1} &= \delta_n + z_n  \\
s_{3,n+1} &= s_{2,n}  \\
&\vdots \\
s_{\nu+1,n+1} &= s_{\nu,n}  
\end{split}
\end{align}
if $\nu>0$.
The slowly oscillating solutions ($\nu=0$) are governed by the map
\begin{align} \label{FreqPresMapk0Compl}
s_{1,n+1} &=  -1 +2 \cos(\omega z_n)  \, e^{-\mu z_n} \nonumber \\
& \quad -  (s_{1,n} + 1) \, \cos(\omega [1 + z_n]) \, e^{-\mu [1+ z_n]}  .
\end{align}

\subsection{Fixed Points}

Symmetric periodic solutions with discrete period $\nu$ are fixed points of $\mathcal{M}_\nu$. Consistent with the symmetry requirement, the simple structure of $\mathcal{M}_\nu$ immediately confirms that the fixed point solutions have equal time-intervals between $x$-zero-crossings. We denote this time interval, the switching interval, by $T^*$. It is half of the period $P_\nu$ of the corresponding four-symbol symmetric periodic solution, i.e. $T^* = P_\nu / 2$.
%
%\begin{align}
% T_{j,n} = T^* = \frac{P_\nu}{2}  \qquad  j=1,\ldots,\nu.
%\end{align}
%
The period satisfies the constraint
\begin{align}
\frac{2}{\nu+1} < P_\nu< \frac{2}{\nu}.
\end{align}
The switching time is the sum of the non-negative time $\delta^*$ to the next sign-switch of the feedback and the  non-negative time $z^*$ to the subsequent zero crossing, $T^* = z^* + \delta^*$.

In the underdamped regime ($Q>1/2$), the ODE-flow $\Phi_\pm$ crosses the switching manifold repeatedly, and one needs the inequalities
\begin{align} \label{zcond}
z^* &= (\nu+1) \,T^* -1 < \pi/\omega
\end{align}
and 
\begin{align} \label{deltacond}
\delta^* = 1- \nu \,T^* < \pi/\omega.
\end{align}
to ensure that the map's fixed point corresponds to a symmetric periodic solution with alternating $H$-type and $Z$-type events, as required by the assumptions made in deriving the map $\mathcal{M}_\nu$. 
 Recalling that $\omega = \Omega \, \sqrt{4 Q^2 -1} /(2 Q)$, the inequalities mean that, for a fixed non-negative integer $\nu$, the corresponding symmetric periodic solution only exists in some region of $\Omega, Q$ parameter space (see Fig.~\ref{fig:RegionPlot} and Fig.~\ref{fig:PeriodPlotNF}a).

In the overdamped regime ($0<Q<1/2$), the ODE-flow $\Phi_\pm$ can cross the switching manifold at most once and no limits on $T^*$ exist. Instead there are limits on the discrete frequency, $\nu$ must be even if the feedback is negative ($\sigma=-1$) and odd if the feedback is positive ($\sigma=+1$).
The implication is that there exists a symmetric periodic 4-symbol solutions for any even (odd) non-negative integer $\nu$ in the case of negative (positive) feedback. The  seven lowest frequency modes for negative feedback are shown in Fig.~\ref{fig:PeriodPlotNF}b.

Assuming above conditions are satisfied, the state-vector of the fixed-point is 
\begin{align} \label{fixedpoint}
\mathbf{s}^* = \left( y_Z^*, T^*, \ldots, T^* \right)^T,
\end{align}
($\mathbf{s}^* =  y_Z^*$ if $\nu=0$)
with the switching interval $T^*$ given implicitly  by the (smallest) positive root of 
\begin{align}\label{Tfixed}
\tan\left( \omega  [(\nu+1) T^*- 1] \right) = \frac{\sin(\omega  T^*)}{e^{\mu T^*} + \cos(\omega  T^*)}
\end{align}
%which holds for all non-negative integers $\nu$,
%which can alternatively be written as the requirement that $\sin(\omega z^*) \exp(\mu z^*) = \sin(\omega \delta^*) \exp(-\mu \delta^*)$.
and the y-coordinate of the $Z$-type event being
\begin{align} \label{yfixed}
y^*_{Z} =-1 + \frac{ 2}{ e^{\mu z^*} \, \cos(\omega z^*) +  e^{- \mu \delta^*} \, \cos(\omega \delta^*) }.
\end{align}

\subsection{Corner Collisions}

In terms of the fixed points of map $\mathcal{M}_\nu$, corner collisions occur if one of the conditions \eq{zcond} and \eq{deltacond} is violated. In such a case, there no longer exists a valid fixed point of $\mathcal{M}_\nu$. The corresponding periodic solution of the \dde{DDEmodel} may cease to exist due to a bifurcation or it may continue to exist but correspond to a fixed point of a different map, such as $\mathcal{M}_{\nu'}$ with $\nu' \ne \nu$. $\mathcal{M}_\nu$ exhibits both types of corner collisions:

(1) As $\omega$ approaches $\omega \to (\nu+1) \pi$ from below, we find that condition~(\ref{deltacond}) is violated because $\delta^* \to \frac{1}{\nu+1}$  (and $z^* \to 0$). The $H/\overline{H}$-event point of the symmetric periodic solution approaches the switching manifold and collides with the $Z/\overline{Z}$-event point of the periodic solution.
That is, as  the parameter $\omega$ is increased from below to above $\omega = (\nu+1) \pi$, the headpoint of the  4-symbol symmetric periodic moves through the switching manifold and the ordering of the 4-symbol sequence changes due to $Z/\overline{Z}$ and $H/\overline{H}$ symbols exchanging places. Furthermore, the number of zero crossings $\nu$ in the unit delay interval increases by one, such that one needs to consider the fixed points of the map $\mathcal{M}_{\nu+1}$ in order to continue the periodic solution. 

(2) As $\omega$ approaches $\omega \to (2 \nu +1) \pi$ from below, we find that conditions~(\ref{zcond}) and (\ref{deltacond}) are violated because $z^* \to \frac{1}{2 \nu+1}$  and $\delta^* \to \frac{1}{2 \nu +1} $.
The $H/\overline{H}$-event point approaches the switching manifold but does not collide with a previously existing $Z/\overline{Z}$-event point.  Instead,  the 4-symbols symmetric periodic solution collides with another ``nearby'' periodic solutions and ceases to exist.

 Nearby periodic solutions must be solutions with symbol sequences of more than 4-symbols per period, such as $8, 12, 16, \ldots$ symbols. We find that there exists a symmetric 12-symbol periodic solution that coexists with the 4-symbol periodic solution and collides with it at the critical value of $\omega$. One may view this 12-symbol solution as a fixed point of the third iterate of $\mathcal{M}_\nu$, one that is distinct from the fixed point of $\mathcal{M}_\nu$. We omit the algebra, but show in Fig.~\ref{fig:BifPlotk2k3} the numerical continuation of a 12-symbol solution that was obtained analytically.

\subsection{Smooth Bifurcations}

In addition to discontinuity induced transitions, the stability of a symmetric 4-symbol periodic solutions can change due to smooth bifurcations. These correspond to standard bifurcations of the fixed points of the map $\mathcal{M}_\nu$.  Bifurcation curves in parameter space are found by determining the characteristic roots of $\mathcal{M}_\nu$ linearized about the fixed point. Since the Poincare map of symmetric periodic solutions is $\mathcal{P}_\nu = \mathcal{M}_\nu^2$, it is the square of a characteristic root of $\mathcal{M}_\nu$ that determines the bifurcation type.

We find that the Jacobian of the map $\mathcal{M}_\nu$, given by \eq{FreqPresMapCompl}, has the form
\begin{align}  \label{Jac}
\mathbf{D\mathcal{M}}_\nu = \begin{pmatrix}
a  &  b   & b &  \ldots & b &b & b  \\
c  &  d   & d & \ldots & d & d & d \\
0 & 1 & 0 &  \ldots & 0 & 0 & 0  \\
0 & 0& 1 &  \ldots & 0 & 0 & 0\\
&&& \ddots &&& \\
0 & 0 & 0 &  \ldots & 1 & 0 & 0  \\
0 & 0 & 0 &  \ldots & 0 & 1 & 0  \\
\end{pmatrix}
\end{align}
with coefficients given in App.~\ref{appendix_Jac_coeff}.
The $(\nu+1)$ roots $\lambda$ that satisfy
\begin{align} \label{CharEqMatrix}
0 =  \left| \mathbf{D\mathcal{M}}_\nu - \lambda \mathbf{I} \right|
\end{align}
are solutions to the characteristic equation
\begin{align} \label{CharEq}
0 = [(a\!- \!1)  d - b  c]  \frac{1-\lambda^\nu}{1-\lambda} +  d - (a + d)  \lambda^\nu + \lambda^{\nu+1},
\end{align}
as shown in App.~\ref{appendix_chareq}.

We first consider $\nu=0$ periodic solutions, slowly oscillating solutions.  In this case, the characteristic equation reduces to $\lambda=a$. It can be shown that $|a|<1$ (see App.~\ref{appendix_a_bound}), implying that the characteristic roots have magnitude less than one independent of the choice of parameters. The slowly oscillating solutions are always locally stable.

We consider next bifurcations of the rapidly oscillating symmetric 4-symbol periodic solutions ($\nu>0$), focusing on those bifurcations that were found numerically, namely pitchfork and Neimark-Sacker bifurcations. We provide expressions that not only allow the bifurcation curves to be determined analytically but enable us to specify parameter regions where bifurcations cannot occur.

Pitchfork bifurcations are associated with a characteristic root equal to $+1$ of the Poincar\'{e} map $\mathcal{P}_\nu$ and arise instead of a generic saddle-node bifurcation due to the inversion-symmetry of \dde{DDEmodel}. One needs to consider, therefore, whether there exist curves in parameter space along which $\mathcal{M}_\nu$ has a characteristic root that is either $\lambda=1$ or $\lambda=-1$. As shown in App.~\ref{appendix_nounitylambda}, there exists no solution $\lambda=1$ of \eq{CharEq}; pitchfork bifurcations are associated with a characteristic value $\lambda=-1$. In this case, \eq{CharEq} reduces to
\begin{subequations} \label{pitchfork}
\begin{align}
 0 &=  -a-1 & \text{if } \nu &\text{ even} \label{pitchfork_even}\\
0 &= (1+d) + e^{-2 \mu T^*} & \text{if } \nu &  \text{ odd} \label{pitchfork_odd},
\end{align}
\end{subequations}
where we made use of the identity \eq{abcd_identity_2}. The equality \eq{pitchfork_even}, covering the case of even $\nu$, cannot be satisfied because $|a|<1$. The equality \eq{pitchfork_odd}, covering the case of odd $\nu$, cannot be satisfied in the overdamped regime nor in the underdamped regime if $\omega T^* < \pi$  because under these conditions it can be shown that  $d+1 > -e^{-2 \mu T^*}$ (see App.~\ref{appendix_d_bound}).  Thus, a pitchfork bifurcation can occur only if the requirements (1) $\nu$ is odd and (2) $\omega T^* > \pi$ are simultaneously satisfied, which is only possible for negative feedback.  No pitchfork bifurcation of symmetric periodic solutions can occur if the feedback is positive. For negative feedback, an example of a pitchfork bifurcation curve found by solving \eq{pitchfork_odd} is shown in Fig.~\ref{fig:BifPlotk2k3}a (the analytically determined curve and bifurcation curve obtained numerically coincide). As seen in this example, the pitchfork bifurcation of the symmetric periodic solution gives rise to a pair of asymmetric periodic solutions (Fig.~\ref{fig:BifPlotk2k3_Asym}).

Neimark-Sacker bifurcations are associated with pairs of complex conjugate roots of magnitude one.  Accordingly, we seek parameter values for which $\lambda=e^{i \phi}$ with $\phi \in (0,\pi)$ is a solution to the characteristic equation. After some algebraic manipulation and separation of imaginary and real parts, we obtain 
\begin{subequations} \label{HopfBif}
\begin{align}
0 &=  \Big[ f_1 + \cos \phi \Big] \frac{ \sin(\frac{\nu+1}{2} \phi) }{\sin(\frac{\phi}{2}) } + f_2 \, \cos(\frac{\nu}{2} \phi) 
 \label{HopfBifCharEqRe} \\
0 &= \sin(\frac{\nu+1}{2}\phi) \cos(\frac{\phi}{2})  +  f_3 \, \sin(\frac{\nu}{2} \phi),  \label{HopfBifCharEqIm}
\end{align}
\end{subequations}
where $f_1, f_2$ and $f_3$ are functions of $\Omega$ and  $Q$ (see App.~\ref{appendix_NS}).
We utilize \eq{HopfBifCharEqIm} to obtain $\phi$ as a function of the two parameters,  $\phi = \phi_n(\Omega, Q)$, where $n$ labels the complex-root pair. Substitution of $\phi_n$ into \eq{HopfBifCharEqRe} allows us then to find the Neimark-Sacker bifurcation curves in the $(\Omega, Q)$ parameter plane for each of the complex-root pairs.

Examples of Neimark-Sacker bifurcation curves are seen in Fig.~\ref{fig:BifPlotk2k3}a for the case of negative feedback.   There is one curve for the $\nu=2$ periodic solution, which has a 3-dimensional Poincar\'{e} map, and two curves for the $\nu=3$ periodic solution,  which has a 4-dimensional Poincar\'{e} map.

In the overdamped regime, symmetric solutions with $\nu>0$ are unstable for small $\Omega$ and become stable after undergoing an appropriate number of Neimark-Sacker bifurcations as $\Omega$ is increased (see Fig.~\ref{fig:PeriodPlotNF}b).
Numerically, we find that the $\nu+1$ roots $\lambda$ of the characteristic equation for $Q<1/2$ remain inside the unit circle in the limit $\Omega \to \infty$.
Thus, in the overdamped regime, each mode first becomes stable and then retains stability as $\Omega$ is increased. The number of stable coexisting solutions grows with $\Omega$. That is, for any chosen discrete frequency $\nu^*$, there is a sufficiently large $\Omega$ such that all symmetric 4-symbol periodic solutions with even (odd) $\nu$ smaller or equal to $\nu^*$ are stable and coexist if the feedback is negative (positive).

\section{Discussion} 

In this paper we advance the understanding of  periodic solutions that arise in second order linear DDEs with relay feedback. The DDE discussed is a representative model of systems with a delayed and bandpass filtered relay-type feedback signal. It also represents applications that exhibit approximate harmonic oscillator type dynamics and have a time delayed relay feedback of the velocity signal. 

We show that it is useful to distinguish the underdamped and overdamped regime. In the overdamped regime all stable solutions found are periodic, whereas a much richer solution and bifurcation structure exists in the underdamped regime.  For example, when underdamped, stable periodic orbits can coexist with stable quasiperiodic solutions. In either regime, the system exhibits strong multirhythmicity. In terms of applications in which such systems serve as signal generators, this means that they are able to produce a large number of distinct periodic modes. These modes can be accessed either by controlling initial conditions or through parameter tuning.

Similar to first order delayed relay systems, the slowly oscillating solution of \dde{DDEmodel} is always stable if it exists. The slowly oscillating solution exists for all parameters in the overdamped regime if the feedback is negative. It also exists in the underdamped regime, both for positive and negative feedback, if $\Omega$ is sufficiently small. 

In this paper we have restricted our investigation to a linear second order DDE with symmetry. This significantly simplified the analytic treatment. If symmetry is lifted, then smooth bifurcations are expected to unfold in the usual way. For example, pitchfork bifurcations of periodic-solution fixed-points will unfold into corresponding saddle-node bifurcations, similar to what is found in \cite{Barton2006}.
The lifting of the symmetry, for example by moving the ODE fixed points off of the switching manifold, would also affect discontinuity induced transitions. Tangential grazing bifurcation~\cite{Sieber2006} become possible~\cite{Benadero2019}. 
The techniques described in this paper can be extended in a straightforward way to the asymmetric case.

The linearity of the ODEs governing the dynamics of \dde{DDEmodel} permits an explicit construction of finite dimensional maps. For nonlinear ODEs, maps can be constructed near periodic orbits~\cite{Sieber2006} but global results are more difficult to obtain.  Nevertheless, the behavior of the linear system is a good starting point for studies of related models with nonlinearities.

An interesting extension of our work would be to consider relays that feature intrinsic hysteretic behavior as well as delay~\cite{Sieber2010}, as this is a good model of many control elements used in practice. 
It would also be fruitful to explore the dynamics of \dde{DDEmodel} with the step-like relay nonlinearity replaced by a smoothed version, because infinitely sharp step functions are not achievable in most applications. 
One often finds good correspondence.
As an example, climate phenomena described by a DDE model containing a sigmoidal type nonlinearity was studied numerically in~\cite{Keane2015,Keane2016} and the numerically observed behavior of the smooth DDE could be explained by analysis of a nonsmooth DDE that resulted from replacing the sigmoidal nonlinearity with a switching function~\cite{Ryan2020}.  
Similarly,  a second order DDE  related to the pupil light reflex was investigated hand-in-hand with a related smooth system in~\cite{Barton2006}. While the dynamics in the vicinity of the nonsmooth bifurcations is strongly changed under transition to the smoothed system, these changes happen in a controlled manner allowing one to establish a clear connection. It would be valuable to explore the correspondence of smooth and nonsmooth dynamics for \dde{DDEmodel}.

\begin{acknowledgments}
This study was supported in part by the Irish Research Council and McAfee LLC under award number EBPPG/2018/269.
 \end{acknowledgments}

\appendix

\section{Coefficients of Jacobian \label{appendix_Jac_coeff}}

The coefficients of the Jacobian matrix of the map $\mathcal{M}_\nu$ are
\begin{align}
a &= - e^{-\mu T^*} \left[ \cos(\omega T^*) + \frac{\mu}{\omega} \sin(\omega T^*) \right] \label{a_coeff}
\\
b &=  - (\mu^2 + \omega^2) (y^*_{Z}  + 1)  \frac{e^{-\mu T^*}}{\omega} \sin(\omega T^*)  \label{b_coeff}
\\
c &=   \frac{1}{(y^*_{Z}  + 1)} \frac{e^{-\mu T^*}}{\omega}  \sin(\omega T^*)   \label{c_coeff}
\\
d&= -1 - \left[   \cos(\omega T^*)  - \frac{\mu}{\omega} \sin(\omega T^*)  \right] \, e^{- \mu T^*} \label{d_coeff}
\end{align}
and we note the following useful identities
\begin{align}
[(a-1) d - b c]  &=  1 + 2 \cos(\omega T^*) \, e^{-\mu T^*} + e^{-2 \mu T^*}  \label{abcd_identity_1}\\
[a (d+1) - b c] &=  e^{-2 \mu T^*}.  \label{abcd_identity_2}
\end{align}

\section{Computation of the characteristic equation \label{appendix_chareq}}

Consider the determinant of the $\nu \times \nu$ matrix
\begin{align*}
A =
 \begin{pmatrix}
 -\lambda \,\alpha    & 0 &  \ldots & 0 &0 & 1   \\
 1 & -\lambda &  \ldots & 0 & 0 & 0  \\
 0& 1 &  \ldots & 0 & 0 & 0\\
&& \ddots &&& \\
 0 & 0 &  \ldots & 1 & -\lambda & 0  \\
 0 & 0 &  \ldots & 0 & 1 & -\lambda  \\
 \end{pmatrix}.
 \end{align*}
If, when evaluating the determinant, one pulls out $\alpha$,  pulls out $\lambda^{-1}$ from the second row, $\lambda^{-2}$ from the third row, $\lambda^{-3}$ from the forth row, etc., and, subsequently, adds the first to the second row, the second to the third row, and so forth; then one obtains 
\begin{align*}
 \det A &=  \frac{\alpha}{\lambda^{(\nu-1) \nu/2} }   \begin{vmatrix}
 -\lambda    & 0 &   \ldots  &0 &\alpha^{-1}    \\
0  & -\lambda^2 &   \ldots  & 0 & \alpha^{-1} \\
&& \ddots &&& \\
 0 & 0 & \ldots  & -\lambda^{\nu-1}  & \alpha^{-1}   \\
 0 & 0 &  \ldots & 0 & \alpha^{-1} -\lambda^\nu  \\
 \end{vmatrix},
\end{align*}
which is upper triangular and evaluates to
\begin{align} \label{SimpleDet}
 \det A &= (-1)^{\nu-1} (1-\alpha \lambda^\nu).
\end{align}
The determinant 
\begin{align*}
\Delta(\beta) &= \begin{vmatrix}
 \beta    & 1 &  \ldots & 1 &1 & 1   \\
 1 & -\lambda &  \ldots & 0 & 0 & 0  \\
 0& 1 &  \ldots & 0 & 0 & 0\\
&& \ddots &&& \\
 0 & 0 &  \ldots & 1 & -\lambda & 0  \\
 0 & 0 &  \ldots & 0 & 1 & -\lambda  \\
 \end{vmatrix}
 \end{align*}
 can be obtained by pulling out  $(1-\lambda)^{-1}$ from the first row, followed by subtraction of all rows other than the first from the first row, yielding
\begin{align*}
\Delta(\beta) =
\frac{1}{1-\lambda} 
 \begin{vmatrix}
 \beta -1-\beta \lambda & 0 & \ldots & 0 & 0 & 1 \\
 1 & -\lambda &  \ldots & 0 & 0 & 0  \\
 0& 1 &  \ldots & 0 & 0 & 0\\
&&& \ddots &&& \\
 0 & 0 &  \ldots & 1 & -\lambda & 0  \\
 0 & 0 &  \ldots & 0 & 1 & -\lambda  \\
 \end{vmatrix}
 \end{align*}
 which, upon utilizing \eq{SimpleDet}, evaluates to
 \begin{align} 
\Delta(\beta) = (-1)^{\nu-1} \left[ \frac{1-\lambda^{\nu-1}}{1-\lambda}  + \beta \lambda^{\nu-1} \right]. \label{SimpleDet2}
\end{align}
Utilizing Laplace expansion on the first row of the matrix in \eq{CharEqMatrix},  for which the Jacobian is given by \eq{Jac}, and then pulling out, respectively, $d$ and $b$, we reduce the problem to 
\begin{align*}
0 =(a-\lambda)  d \; \Delta(1-\lambda/d) - c \, b \; \Delta(1),
 \end{align*}
which  yields the characteristic equation
\begin{align*}
0 =[(a-1) \, d - b \, c] \; \frac{1-\lambda^\nu}{1-\lambda} +  d - (a + d) \, \lambda^\nu + \lambda^{\nu+1}.
\end{align*}

\section{Bound on $a$ \label{appendix_a_bound}}

To show that $|a|<1$, we treat the underdamped and overdamped case separately: 
In the underdamped case both $\mu$ and $\omega$ are positive real and it is  seen that $|a| <1$ because
\begin{align}
|a| &= \left|e^{-\mu T^*} \right| \left|\cos(\omega T^*) + \frac{\mu}{\omega} \sin(\omega T^*) \right| \nonumber \\
&\le \left|e^{-\mu T^*} \right| \left(  \left|\cos(\omega T^*) \right| + (\mu T^*) \left| \frac{\sin(\omega T^*)}{\omega T^*} \right| \right) \nonumber \\
& \le e^{-\mu T^*}  \left(  1 + \mu T^*  \right) \nonumber \\
&< 1.
\end{align}
If $Q<1/2$, such that the system is overdamped, $\omega$ is imaginary ($\omega = i |\omega|$) and $\mu > |\omega|$ holds,  as seen from \eq{muomdef}. Furthermore, $a<0$, such that $|a| = -a$. We can write
\begin{align*}
|a|  =& \, e^{-\mu T^*} \left[ \cosh(|\omega| T^*) +  (\mu T^*)  \frac{\sinh(|\omega| T^*)}{|\omega| T^*} \right] \\
=& \, e^{-(\mu - |\omega|) T^*} \left[ 1 + (\mu - |\omega|) T^* \right] \nonumber \\ 
&- \frac{(\mu - |\omega|) \, e^{-(\mu - |\omega|) T^*} }{2 |\omega|} \left[ e^{- 2 |\omega| T^*} + 2 |\omega| T^* - 1\right]. 
\end{align*}
Since $0 \le  e^{- 2 |\omega| T^*} + 2 |\omega| T^* - 1$ with equality for $ |\omega| =0$, we find
\begin{align}
|a| \le e^{-(\mu - |\omega|) T^*} \left[ 1 + (\mu - |\omega|) T^* \right] < 1.
\end{align}

\section{No $\lambda=1$ solution \label{appendix_nounitylambda}}

There exists no solution $\lambda=1$ because substitution of $\lambda=1$ into \eq{CharEq} produces the condition
\begin{align}
0 =[(a-1) \, d - b \, c] \; \nu + \left[1 -a \right],
\end{align}
which cannot be satisfied because it has a positive right hand side. The first term on the right hand side is positive because, after rewriting $[(a-1) \, d - b \, c]$ by using \eq{abcd_identity_1}, one finds
\begin{align}
  1 + 2 \cos(\omega T^*) \, e^{-\mu T^*} + e^{-2 \mu T^*} \ge (1-  e^{-\mu T^*})^2.
\end{align}
The second term is positive because $|a|<1$ as shown in App.~\ref{appendix_a_bound}.

\section{Bound on $d$ \label{appendix_d_bound}}

The purpose of this section is to show that if either $Q\ge1/2$ and $\omega T^* < \pi$ or $Q<1/2$ then
\begin{align}
-e^{-2 \mu T^*} < d+1.
\end{align}
To demonstrate this inequality we consider the expression 
\begin{align} \label{dplus1a}
g = e^{\mu T^*} (d+1)&=  \frac{\mu}{\omega} \sin(\omega T^*) -   \cos(\omega T^*),
\end{align}
where we made use of the definition of $d$ given by \eq{d_coeff}. 
We will show that $g$ is bounded from below by 
\begin{align}
f = - e^{-\mu T^*}.
\end{align}
We treat the overdamped and underdamped regime separately:
In the overdamped regime ($Q<1/2$), $\omega = i |\omega| = i \mu \sqrt{1-4 Q^2}$.  For convenience, we introduce $x = |\omega| T^*$  and $y =  {\mu}/{|\omega|}$ in order to rewrite \eq{dplus1a} as
\begin{align}
g &=  y \sinh(x) -   \cosh(x)  
\end{align}
with $x \in (0,\infty)$ and $y \in (1,\infty)$ and $f$ as $f= - e^{-x y}$.
First note that $g$ and $f$ are equal at $x=0$ and that both $f$ and $g$ have positive slope with respect to $x$ with $g' > f'$ because $g' = (y-1) \sinh(x) + y e^{-x}$, $f' =  y  e^{-x y}$, and 
\begin{align}
(y-1) \sinh(x) + y e^{-x} > y e^{-x}> y e^{-x y}  >0
\end{align}
if $y>1$ and $x>0$. Thus, $g>f$ if $x>0$ and $y>1$.

If $Q\ge1/2$, $\omega$ is real and we introduce the variable $\phi = \omega T^*$ together with the assumption that $0\le\phi<\pi$. For convenience we, furthermore, introduce $\mathtt{y}=\mu T^*$ and write \eq{dplus1a} as
\begin{align}
g &=  \mathtt{y} \, \mathrm{sinc}(\phi) -   \cos(\phi)  \qquad \phi \in [0,\pi), \, \mathtt{y} \in (0,\infty).
\end{align}
Utilizing the inequalities $e^{-\mathtt{y}} > (1-\mathtt{y})$ if $\mathtt{y}\ne0$ and $1 \ge \mathrm{sinc}(\phi)>0$ if $ \phi \in [0,\pi)$, yields the estimate
\begin{align}
-f = e^{-\mathtt{y}} > (1 - \mathtt{y}) \ge (1 - \mathtt{y}) \, \mathrm{sinc}(\phi).
\end{align}
The bound on $d$ is obtained by noting that $\mathrm{sinc}(\phi) \ge \cos(\phi)$ if $\phi \in [0,\pi)$, such that 
\begin{align}
-f >   \cos(\phi) - \mathtt{y} \, \mathrm{sinc}(\phi) = - g.
\end{align}

\section{Neimark-Sacker bifurcation \label{appendix_NS}}

Here we derive \eq{HopfBif}. First note that the characteristic equation, \eq{CharEq}, can be rewritten as
 \begin{align*}
0 = A \frac{\lambda^{\nu+1}-1}{\lambda-1}  +(\lambda^{\nu+1} - 1 )
- B  (\lambda^\nu-1)
 +  d +1  - B 
 \end{align*}
 with $A =[(a-1) d - b \, c]$ and $B = [(a (d+1) - b c]$. Substitution of $\lambda=e^{i \phi}$, multiplication by $e^{-i \nu \phi/2}$, and subsequent separation of real and imaginary parts gives the equation pair
  \begin{align*}
0 &= [f_1 + \cos \phi ] \;  \frac{ \sin(\frac{\nu+1}{2} \phi) }{\sin(\frac{\phi}{2})}  + f_2 \cos(\frac{\nu}{2} \phi) \\
0 &= \cos(\frac{\phi}{2})\,\sin(\frac{\nu+1}{2} \phi) + f_3 \sin(\frac{\nu}{2} \phi).
 \end{align*}
For notational convenience, we introduced $f_1 = A - 1$, $f_2 = d+1-B$, and $f_3= -d-1-B$. Utilizing the definition~(\ref{d_coeff}) of $d$ as well as the identities \eq{abcd_identity_1} and \eq{abcd_identity_2} and recalling that $\omega$ and $\mu$ are defined in terms of the parameters $\Omega$ and $Q$ according to \eq{muomdef}, we find that the function $f_i(\Omega, Q)$ are given by
\begin{align}
f_1 &=  2 \cos(\omega T^*) \, e^{-\mu T^*} + e^{-2 \mu T^*}\\
%&
f_2 &=  - \left[   \cos(\omega T^*)  - \frac{\mu}{\omega} \sin(\omega T^*)  \right] \, e^{- \mu T^*} - e^{- 2 \mu T^*} \\
f_3 &=   \left[   \cos(\omega T^*)  - \frac{\mu}{\omega} \sin(\omega T^*)  \right] \, e^{- \mu T^*} - e^{- 2 \mu T^*} .
\end{align}

%%%
\bibliography{NonSmoothDelaySystem}

\begin{thebibliography}{16}
\expandafter\ifx\csname natexlab\endcsname\relax\def\natexlab#1{#1}\fi
\expandafter\ifx\csname bibnamefont\endcsname\relax
  \def\bibnamefont#1{#1}\fi
\expandafter\ifx\csname bibfnamefont\endcsname\relax
  \def\bibfnamefont#1{#1}\fi
\expandafter\ifx\csname citenamefont\endcsname\relax
  \def\citenamefont#1{#1}\fi
\expandafter\ifx\csname url\endcsname\relax
  \def\url#1{\texttt{#1}}\fi
\expandafter\ifx\csname urlprefix\endcsname\relax\def\urlprefix{URL }\fi
\providecommand{\bibinfo}[2]{#2}
\providecommand{\eprint}[2][]{\url{#2}}

\bibitem[{\citenamefont{Barton et~al.}(2006)\citenamefont{Barton, Krauskopf,
  and Wilson}}]{Barton2006}
\bibinfo{author}{\bibfnamefont{D.~A.~W.} \bibnamefont{Barton}},
  \bibinfo{author}{\bibfnamefont{B.}~\bibnamefont{Krauskopf}},
  \bibnamefont{and} \bibinfo{author}{\bibfnamefont{R.~E.}
  \bibnamefont{Wilson}}, \bibinfo{journal}{Dyn. Syst.}
  \textbf{\bibinfo{volume}{21}}, \bibinfo{pages}{289} (\bibinfo{year}{2006}).

\bibitem[{\citenamefont{Ryan et~al.}(2020)\citenamefont{Ryan, Keane, and
  Amann}}]{Ryan2020}
\bibinfo{author}{\bibfnamefont{P.}~\bibnamefont{Ryan}},
  \bibinfo{author}{\bibfnamefont{A.}~\bibnamefont{Keane}}, \bibnamefont{and}
  \bibinfo{author}{\bibfnamefont{A.}~\bibnamefont{Amann}},
  \bibinfo{journal}{Chaos} \textbf{\bibinfo{volume}{30}},
  \bibinfo{pages}{023121} (\bibinfo{year}{2020}), ISSN
  \bibinfo{issn}{1054-1500}.

\bibitem[{\citenamefont{Fridman et~al.}(2000)\citenamefont{Fridman, Fridman,
  and Shustin}}]{Fridman2000}
\bibinfo{author}{\bibfnamefont{E.}~\bibnamefont{Fridman}},
  \bibinfo{author}{\bibfnamefont{L.}~\bibnamefont{Fridman}}, \bibnamefont{and}
  \bibinfo{author}{\bibfnamefont{E.}~\bibnamefont{Shustin}},
  \bibinfo{journal}{J. Dyn. Syst. Meas. Control}
  \textbf{\bibinfo{volume}{122}}, \bibinfo{pages}{732} (\bibinfo{year}{2000}).

\bibitem[{\citenamefont{Fridman et~al.}(2002)\citenamefont{Fridman, Fridman,
  and Shustin}}]{Fridman2002}
\bibinfo{author}{\bibfnamefont{E.}~\bibnamefont{Fridman}},
  \bibinfo{author}{\bibfnamefont{L.}~\bibnamefont{Fridman}}, \bibnamefont{and}
  \bibinfo{author}{\bibfnamefont{E.}~\bibnamefont{Shustin}},
  \emph{\bibinfo{title}{Sliding Mode Control in Engineering}}
  (\bibinfo{publisher}{Marcel Dekker, New York}, \bibinfo{year}{2002}), chap.
  \bibinfo{chapter}{Steady modes in relay systems with delay}, pp.
  \bibinfo{pages}{263--293}.

\bibitem[{\citenamefont{Shustin et~al.}(2003)\citenamefont{Shustin, Fridman,
  and Fridman}}]{Shustin2003}
\bibinfo{author}{\bibfnamefont{E.}~\bibnamefont{Shustin}},
  \bibinfo{author}{\bibfnamefont{E.}~\bibnamefont{Fridman}}, \bibnamefont{and}
  \bibinfo{author}{\bibfnamefont{L.}~\bibnamefont{Fridman}},
  \bibinfo{journal}{Discrete Contin. Dyn. Syst.} \textbf{\bibinfo{volume}{9}},
  \bibinfo{pages}{339} (\bibinfo{year}{2003}).

\bibitem[{\citenamefont{Benadero et~al.}(2019)\citenamefont{Benadero, Aroudi,
  and Ponce}}]{Benadero2019}
\bibinfo{author}{\bibfnamefont{L.}~\bibnamefont{Benadero}},
  \bibinfo{author}{\bibfnamefont{A.~E.} \bibnamefont{Aroudi}},
  \bibnamefont{and} \bibinfo{author}{\bibfnamefont{E.}~\bibnamefont{Ponce}},
  \bibinfo{journal}{Nonlinear Theory Appl. IEICE}
  \textbf{\bibinfo{volume}{10}}, \bibinfo{pages}{337} (\bibinfo{year}{2019}).

\bibitem[{\citenamefont{Sieber}(2006)}]{Sieber2006}
\bibinfo{author}{\bibfnamefont{J.}~\bibnamefont{Sieber}},
  \bibinfo{journal}{Nonlinearity} \textbf{\bibinfo{volume}{19}},
  \bibinfo{pages}{2489} (\bibinfo{year}{2006}).

\bibitem[{\citenamefont{Hale et~al.}(2002)\citenamefont{Hale, Magalh\~{a}es,
  and Oliva}}]{Hale2002}
\bibinfo{author}{\bibfnamefont{J.}~\bibnamefont{Hale}},
  \bibinfo{author}{\bibfnamefont{L.~T.} \bibnamefont{Magalh\~{a}es}},
  \bibnamefont{and} \bibinfo{author}{\bibfnamefont{W.~L.} \bibnamefont{Oliva}},
  \emph{\bibinfo{title}{Dynamics in Infinite Dimensions}},
  vol.~\bibinfo{volume}{47} of \emph{\bibinfo{series}{Applied mathematical
  sciences}} (\bibinfo{publisher}{Springer-Verlag, New York},
  \bibinfo{year}{2002}), \bibinfo{edition}{2nd} ed.

\bibitem[{\citenamefont{Weicker et~al.}(2015)\citenamefont{Weicker, Erneux,
  Rosin, and Gauthier}}]{Weicker2015}
\bibinfo{author}{\bibfnamefont{L.}~\bibnamefont{Weicker}},
  \bibinfo{author}{\bibfnamefont{T.}~\bibnamefont{Erneux}},
  \bibinfo{author}{\bibfnamefont{D.~P.} \bibnamefont{Rosin}}, \bibnamefont{and}
  \bibinfo{author}{\bibfnamefont{D.~J.} \bibnamefont{Gauthier}},
  \bibinfo{journal}{Phys. Rev. E} \textbf{\bibinfo{volume}{91}},
  \bibinfo{pages}{012910} (\bibinfo{year}{2015}).

\bibitem[{\citenamefont{an~der Heiden et~al.}(1990)\citenamefont{an~der Heiden,
  Longtin, Mackey, Milton, and Scholl}}]{Heiden1990}
\bibinfo{author}{\bibfnamefont{U.}~\bibnamefont{an~der Heiden}},
  \bibinfo{author}{\bibfnamefont{A.}~\bibnamefont{Longtin}},
  \bibinfo{author}{\bibfnamefont{M.~C.} \bibnamefont{Mackey}},
  \bibinfo{author}{\bibfnamefont{J.~G.} \bibnamefont{Milton}},
  \bibnamefont{and} \bibinfo{author}{\bibfnamefont{R.}~\bibnamefont{Scholl}},
  \bibinfo{journal}{J. Dynam. Differential Equations}
  \textbf{\bibinfo{volume}{2}}, \bibinfo{pages}{423} (\bibinfo{year}{1990}).

\bibitem[{\citenamefont{Milton and Longtin}(1990)}]{Milton1990}
\bibinfo{author}{\bibfnamefont{J.~G.} \bibnamefont{Milton}} \bibnamefont{and}
  \bibinfo{author}{\bibfnamefont{A.}~\bibnamefont{Longtin}},
  \bibinfo{journal}{Vision Research} \textbf{\bibinfo{volume}{30}},
  \bibinfo{pages}{515} (\bibinfo{year}{1990}).

\bibitem[{\citenamefont{Dimopoulos}(2012)}]{Dimopoulos2012}
\bibinfo{author}{\bibfnamefont{H.~G.} \bibnamefont{Dimopoulos}},
  \emph{\bibinfo{title}{Analog Electronic Filters}}
  (\bibinfo{publisher}{Springer}, \bibinfo{year}{2012}),
  chap.~\bibinfo{chapter}{10}, pp. \bibinfo{pages}{397--399}.

\bibitem[{\citenamefont{Govaerts et~al.}(2008)\citenamefont{Govaerts,
  Kuznetsov, Khoshsiar~Ghaziani, and Meijer}}]{Govaerts2008}
\bibinfo{author}{\bibfnamefont{W.}~\bibnamefont{Govaerts}},
  \bibinfo{author}{\bibfnamefont{Y.~A.} \bibnamefont{Kuznetsov}},
  \bibinfo{author}{\bibfnamefont{R.}~\bibnamefont{Khoshsiar~Ghaziani}},
  \bibnamefont{and} \bibinfo{author}{\bibfnamefont{H.~G.~E.}
  \bibnamefont{Meijer}} (\bibinfo{year}{2008}),
  \urlprefix\url{http://sourceforge.net/projects/matcont}.

\bibitem[{\citenamefont{Sieber et~al.}(2010)\citenamefont{Sieber, Kowalczyk,
  Hogan, and Di~Bernardo}}]{Sieber2010}
\bibinfo{author}{\bibfnamefont{J.}~\bibnamefont{Sieber}},
  \bibinfo{author}{\bibfnamefont{P.}~\bibnamefont{Kowalczyk}},
  \bibinfo{author}{\bibfnamefont{S.}~\bibnamefont{Hogan}}, \bibnamefont{and}
  \bibinfo{author}{\bibfnamefont{M.}~\bibnamefont{Di~Bernardo}},
  \bibinfo{journal}{J. Vib. Control} \textbf{\bibinfo{volume}{16}},
  \bibinfo{pages}{1111} (\bibinfo{year}{2010}).

\bibitem[{\citenamefont{Keane et~al.}(2015)\citenamefont{Keane, Krauskopf, and
  Postlethwaite}}]{Keane2015}
\bibinfo{author}{\bibfnamefont{A.}~\bibnamefont{Keane}},
  \bibinfo{author}{\bibfnamefont{B.}~\bibnamefont{Krauskopf}},
  \bibnamefont{and}
  \bibinfo{author}{\bibfnamefont{C.}~\bibnamefont{Postlethwaite}},
  \bibinfo{journal}{SIAM J. Appl. Dyn. Syst.} \textbf{\bibinfo{volume}{14}},
  \bibinfo{pages}{1229} (\bibinfo{year}{2015}).

\bibitem[{\citenamefont{Keane et~al.}(2016)\citenamefont{Keane, Krauskopf, and
  Postlethwaite}}]{Keane2016}
\bibinfo{author}{\bibfnamefont{A.}~\bibnamefont{Keane}},
  \bibinfo{author}{\bibfnamefont{B.}~\bibnamefont{Krauskopf}},
  \bibnamefont{and}
  \bibinfo{author}{\bibfnamefont{C.}~\bibnamefont{Postlethwaite}},
  \bibinfo{journal}{SIAM J. Appl. Dyn. Syst.} \textbf{\bibinfo{volume}{15}},
  \bibinfo{pages}{1656} (\bibinfo{year}{2016}).

\end{thebibliography}
%\bibliographystyle{plain}
%\nocite{*}

%\begin{thebibliography}{29}
%\end{thebibliography}

\end{document}